\documentclass[times,review]{elsarticle} 

\usepackage{amssymb}
\usepackage{amsthm}
\usepackage{verbatim,url}
\usepackage{bm}
\usepackage{amsmath}
\usepackage{color}
\usepackage{algorithmicx,algpseudocode,algorithm}
\usepackage{tikz}
\usepackage{environ}

\usetikzlibrary{graphs,positioning}
\usetikzlibrary{shapes,arrows}
\usetikzlibrary{hobby,decorations.markings}

\usepackage{color}
\usepackage{graphicx}

\usepackage{caption}
\usepackage{subfig}
\usepackage{tikz}
\usetikzlibrary{arrows}

\usepackage{stmaryrd}

\makeatletter
\newsavebox{\measure@tikzpicture}
\NewEnviron{scaletikzpicturetowidth}[1]{%
  \def\tikz@width{#1}%
  \def\tikzscale{1}\begin{lrbox}{\measure@tikzpicture}%
  \BODY
  \end{lrbox}%
  \pgfmathparse{#1/\wd\measure@tikzpicture}%
  \edef\tikzscale{\pgfmathresult}%
  \BODY
}
\makeatother
\makeatletter
\newcommand{\algmargin}{\the\ALG@thistlm}
\makeatother
\newlength{\whilewidth}
\algdef{SE}[parWHILE]{parWhile}{EndparWhile}[1]
  {\parbox[t]{\dimexpr\linewidth-\algmargin}{%
     \hangindent\whilewidth\strut\algorithmicwhile\ #1\ \algorithmicdo\strut}}{\algorithmicend\ \algorithmicwhile}%
\algnewcommand{\parState}[1]{\State%
  \parbox[t]{\dimexpr\linewidth-\algmargin}{\strut #1\strut}}

\newcommand\norm[1]{\left\lVert#1\right\rVert}
\newcommand{\vsh}{v^{\rm sh}}

\journal{}

\begin{document}
\tikzstyle{block} = [rectangle, draw, 
    text width=4.5em, text centered,rounded corners, minimum height=3em, node 
	distance=2cm]
\tikzstyle{line} = [draw, -latex',>=stealth]
\tikzstyle{cloud} = [draw, rectangle,fill=gray!90, text centered,text=black,rounded corners,
		text width=4.5em,
    minimum height=3em]
\tikzstyle{bubble} = [draw, rectangle,fill=gray!20, text centered,rounded corners,
		text width=4.5em,
    minimum height=3em]
\begin{frontmatter}



\title{Variational optimization and data assimilation in chaotic time-delayed systems with automatic-differentiated shadowing sensitivity}

%
\author[label1]{Nisha Chandramoorthy}
\author[label3]{Luca Magri}
\author[label1]{Qiqi Wang}


%
\address[label1]{Massachusetts Institute of Technology,   Center for Computational Science and Engineering,\\
77 Massachusetts Avenue
Cambridge, Massachusetts,02139, USA}
\address[label3]{University of Cambridge, Engineering Department,\\ Trumpington Street, CB2 1PZ, Cambridge, UK}
%
%
%
\begin{abstract}
In this computational paper, we perform sensitivity analysis of long-time (or ensemble) averages in  chaotic regime using the shadowing algorithm. We introduce automatic differentiation to eliminate the tangent/adjoint equation solvers used in the shadowing algorithm. In a gradient-based optimization, we use the computed shadowing sensitivity to minimize different long-time averaged functionals of a chaotic time delayed system by optimal parameter selection. In combined state and parameter estimation for data assimilation, we use the computed sensitivity to predict the optimal trajectory given information from a model and data from measurements beyond the predictability time. The algorithms are applied to a thermoacoustic model. Because the computational framework is rather general, the techniques presented in this paper may be used for sensitivity analysis of ensemble averages, parameter optimization and data assimilation of other chaotic problems, where shadowing methods are applicable. 
\end{abstract}
\end{frontmatter}

\section{Introduction}

Sensitivities are quantitative
measures of the response of model outputs to infinitesimal
changes in inputs, which are crucial to engineering design \cite{sensder1,sensder2}. They are employed
 in parameter
estimation and model selection \cite{sensder8,sensder9}, uncertainty quantification \cite{sensder3,sensder4}, data assimilation \cite{sensder5,sensder6,sensder7}
and design and optimization (see \cite{sensder10,sensder11,sensder12,Magri2019_amr} for 
recent reviews of the applications of sensitivity derivatives in different
engineering disciplines). With the growing ability
to simulate high-dimensional complex dynamics, much research
effort has been invested into commensurately improving sensitivity
analysis methods. Adjoints of mathematical models have been
developed and used successfully in many fields; for example,
adjoint sensitivity analysis in meteorology \cite{weather1, weather2},
aircraft design \cite{sensana1}, systems biology \cite{sensana2}, chemical kinetics \cite{sensder5}, thermo-fluids~\cite{Magri2019_amr}, among others.
The adjoint method is generally used when the input parameter 
space is high-dimensional. Tangent Linear Models (TLMs) and finite
difference methods are also used for sensitivity analysis
when the dimension of the input parameter space is small enough
that the cost of simulating the original dynamics repeatedly is not prohibitive.
In many of these applications, automatic differentiation (AD) has successfully replaced
the TLM or adjoint computations \cite{ad1,ad2}. For example, MITGCM \cite{mitgcm}, a popular
climate model uses \verb+OpenAD+ \cite{mitgcm1,mitgcm2}, an 
open source source-transformation AD software \cite{openad,openAD-combinatorics,openAD1,paul-AD}; and     
Tapenade~\cite{tapenade1,tapenade2} has replaced adjoint 
differentiation in a few industrial-size numerical codes.

In recent times, simultaneous advances in simulation capabilities and
computing power have led to a proliferation of scale-resolving
simulations of chaotic systems \cite{patrickles,des1,des2,des3}. For many of the above target applications, the relevant 
observables, or outputs, in chaotic systems are 
statistically stationary or infinitely long-time averaged 
functions \cite{des1,des2,des3,Huhn2020_jfm,Huhn2020_nody}. Useful gradients of ensemble averages, 
which are equal to infinite time averages in ergodic systems, cannot be obtained by time-averaging 
the instantaneous gradients in chaotic systems \cite{qiqi}. Indeed,
in the infinite time limit, the time averages of the
instantaneous
gradients diverge despite that the ensemble averages of the functions may have
a well-defined gradient. This is because the tangent space of a chaotic attractor is exponentially unstable. Likewise, the sensitivities computed on time-integrating
the adjoint sensitivities diverge exponentially. For the same reasons, other methods, such as TLMs and AD, 
also fail to compute meaningful sensitivities in chaotic systems.
Due to these challenges, sensitivity
analysis of  chaotic systems has not  
developed as much as sensitivity analysis of non-chaotic systems has
\cite[e.g.,][]{angxiu1,angxiu2,qiqi}. 

One approach for the computation of sensitivities of long-time-averaged functionals in chaotic systems is the Least Squares Shadowing \cite{qiqi,shadowing} (LSS) method.  This method  by-passes the exponential instability by computing
the derivatives along a close 
\emph{shadowing} direction, which  is obtained  as the solution of a 
constrained minimization problem. A recent variant 
of the LSS method--the Non-Intrusive Least Squares Shadowing (NILSS) \cite{angxiu1}--
has been proposed to reduce  the computational cost 
and memory requirements of the original LSS problem by projecting the gradients
onto the \emph{unstable} subspaces only.
The NILSS method has found successful applications
in chaotic computational fluid dynamics; e.g.,  \cite{patrickles} and
\cite{angxiu2} applied  NILSS to scale-resolving
Direct Simulations of chaotic
 flows around bluff bodies,  \cite{patrick} developed the adjoint version of NILSS and applied it to wall-bounded
chaotic flows; and \cite{Huhn2020_jfm} applied it to the optimization of chaotic acoustic oscillations subject to synthetic turbulence.  
Other methods for sensitivity analysis
in chaotic systems 
are conceptually based on extensions of the fluctuation-dissipation
theorem for nonequilibrium systems in physics. One  method is based on estimating the 
invariant probability distribution \cite{cooper}. Other
recent approaches \cite{abramov1,colangeli,lucarini,majda} 
 computationally evaluate Ruelle's response formula
for nonequilibrium systems \cite{ruelle2}. In this paper, we are concerned
with the NILSS algorithm and, in particular, on the development
of the automatically
differentiated version of the algorithm. The algorithm is generalized to tackle time delayed systems. 

Delayed differential systems, which often tend 
to be chaotic, are extensively used  
for mathematical modelling of  transport and non-Markovian processes, such as population dynamics
and cell proliferation in mathematical ecology and biology \cite{ddeapp1,ddeapp2}, chemical processes \cite{ddeapp4}, 
neural networks, networked control systems \cite{ddeapp5,ddeapp6} 
thermoacoustics~\cite{Magri2019_amr}, among others. Sensitivity derivatives
with respect to the parameters,
including the delay parameter, have been employed for model selection,
system identification and stability analysis \cite{baker,ddeapp7,Magri2019_amr}. The
objective function to be optimized, in the case of the model parameter estimation or
system identification problems, is typically a mean quantity that
depends on the parameters, including the time delays \cite{baker,Huhn2020_jfm}. The tangent and adjoint discrete AD-shadowing methods developed in this  paper offer
a solution to this problem in chaotic
time-delayed systems.
In order to compute sensitivities using AD, differentiating
the numerical solution of the primal dynamics is necessary. 
If the time delay is treated approximately as an integer 
number of timesteps by the time integrator, differentiating
with respect to the delay poses a     
problem since chain rule differentiability is lost. %
In order to circumvent this issue, in this paper we exploit the
fact that the effect of a  time delay can be replaced
by a linear advection equation, and, therefore, we solve for an extended
primal system. 

A practical engineering problem modelled with time-delayed equations is thermoacoustics~\cite{Magri2019_amr}. 
Gas-turbine and rocket-motor manufacturers strive to design engines that do not experience thermoacoustic instabilities  \cite{Lieuwen2005}. Thermoacoustic instabilities occur when the heat released by the flame is sufficiently in phase with the acoustic pressure~\cite{Rayleigh1878}, such that the thermal energy of the flame that is converted into acoustic energy exceeds dissipation mechanisms.
Unstable thermoacoustic systems have intricate nonlinear behaviours when design parameters are varied, from periodic, through quasi periodic to chaotic oscillations~\cite{Kabiraj2011}. Although methods to investigate the sensitivity of fixed points (with eigenvalue analysis) and periodic solutions (with Floquet analysis) are well-established~\cite{Magri2019_amr}, a stability and sensitivity framework to tackle chaotic acoustic oscillations is only at its infancy~\cite{Huhn2020_jfm,Huhn2020_nody}. 
%
%
%
In themoacoustics, sensitivity analysis quantitatively informs the practitioner on   how to optimally change design parameters, such as geometric quantities;  which passive device is most stabilizing;
and  how large is the uncertainty of the stability calculations~\cite{Magri2016c,Silva2016,Mensah2018}, as reviewed by~\cite{Magri2019_amr}. All these studies are concerned with the calculation of sensitivities of eigenvalues around non-chaotic attractors.  
These established eigenvalue-sensitivity methods fail in chaotic systems because of the butterfly effect~\cite{nisha-ens, lea, Eyink2004,Huhn2020_jfm} (\S\ref{sec:background}). In this paper, we apply the computational framework we develop to the calculation of the derivative of two infinite time-averaged cost functionals, one being an energy norm and the second being an integral metric, with respect to the parameters' vector. These  derivatives give us a quantitative estimate of the long-term response of chaotic acoustic oscillations. We use these sensitivities to stabilize a nonlinearly unstable, yet eigenvalue-stable, thermoacoustic system. Physically, the cost functionals represent the acoustic energy, which we want to minimize to make the combustor operate in stable conditions. We use these sensitivities in a gradient-based optimization algorithm to suppress a chaotic acoustic oscillation, which cannot be achieved by only stabilizing the eigenvalues or through short-term chaotic sensitivity calculations. 

The paper is structured as follows.
In section
\ref{sec:nilssidea}, the idea behind the NILSS algorithm is reviewed.
Section \ref{sec:background}  defines the problem  with a  mathematical background on sensitivity analysis in chaotic systems. The main features of the shadowing algorithm are explained in \ref{sec:nilss}. The AD version of the algorithm is provided in \ref{sec:appxAD}.
We present the chaotic time-delayed model of a prototypical thermoacoustic system in section \ref{sec:rijke}. The tangent and adjoint shadowing sensitivities of this model are calculated and applied for parameter estimation for gradient-based optimization in section \ref{sec:optim}, and for data assimilation in section \ref{sec:da}. The papers ends with a final discussion in Section~\ref{sec:conclusions}. 

\section{Shadowing sensitivity in chaotic systems}
\label{sec:background}
Before we describe the NILSS algorithm, we recall the problem
of extreme sensitivity to perturbations in chaotic systems, which 
leads to ill-conditioning of linearized  models, such as
the tangent equation, the adjoint equation and 
algorithmic differentiation. We define the primal problem 
by a set of ordinary differential equations (ODEs), which 
may be spatially discretized partial differential 
equations, as,
\begin{align}
		\notag
		\dfrac{du}{dt} = \mathcal{F}(u,\mathcal{S}), \; \mathcal{S} \in \mathbb{R}^p \\
		u(0) = u_0 \in \mathbb{R}^d.
		\label{eqn:ode}
\end{align}
Here $u \in \mathbb{R}^d$ is the state of the system, and, 
$\mathcal{S} \in \mathbb{R}^p$ is a vector of system  parameters. The system parameters, which can be, e.g., control variables in an adjoint-based design problem, do not change with time. The 
right hand side $\mathcal{F}:\mathbb{R}^d\times 
\mathbb{R}^p \to \mathbb{R}^d$ of the primal ODE (Eq. \ref{eqn:ode}), which is also referred to as the \emph{time-derivative} direction, is a function of the instantaneous state and $\mathcal{S}$.
In this paper, we study the discrete-time system obtained by time-integration of the 
primal ODE. Throughout, we use a subscript to denote a discrete time, which is 
represented by a positive integer. In particular, $u_0 \in \mathbb{R}^d$ is 
the initial state; $u_n \in \mathbb{R}^d$ is the solution vector at time $n \in \mathbb{Z}^+$. 

We define the function $f:\mathbb{R}^d \times \mathbb{R}^p 
\to \mathbb{R}^d$ to denote the time-one map, i.e., the
time-integrator that evolves a solution state by one timestep, 
so that $u_1 = f(u_0, \mathcal{S}).$
We use the notation $f_n$ to denote the $n$-time 
composition of the map $f,$ at a fixed set of parameters, so that 
$u_n = f_n(u,\mathcal{S})$, $n\in \mathbb{Z}^+$. The set of vectors $\left\{ u_n\right\}$ is an \emph{orbit} or a \emph{trajectory} of the dynamics $f(\cdot, \mathcal{S})$. 
Let $\mathcal{J}$ be a set of $l$ scalar observables in $\mathcal{C}^2(\mathbb{R}^d)$, and $J$ be an observable in this set. Given an initial
state $u_0$, the $N$-time average of $J$ is $\langle J\rangle_N := (1/N)\sum_{n=0}^{N-1} J(u_n).$ In ergodic systems, in the limit $N\to \infty,$ the $N$-time-average, which is referred to as \emph{ergodic average} and denoted as $\langle J\rangle$, is well-defined and independent of the initial state $u_0$. 
The ergodic average $\langle J\rangle$ is a function of the parameters $\mathcal{S}$ only. Its value is equal to an expectation of $J$ with respect to the ergodic, stationary probability distribution achieved by the state vector under the dynamics $f$. In chaotic systems, ergodic averages of observables are often the quantities of interest for  optimization and control problems. In these problems, 
the long-term response of a chaotic system to infinitesimal perturbations may be desired (e.g. \cite{lucarini-longterm}), as opposed to a short-term or intermediate-term response. The problem of nonlinear acoustic oscillations that is studied in this paper is one such example \cite{Huhn2020_jfm}. 
Our goal is to compute, for all $J \in \mathcal{J}$ and all $s \in \mathcal{S}$ the quantity,
\begin{align}
    \label{eqn:qoi}
    d_s \langle J\rangle := d_s \Big(\lim_{N\to \infty} \langle J\rangle_N\Big),  
\end{align}
where $d_s:=d/ds$ denotes the 
differentiation operator with respect to $s$. 
We assume that the ergodic average $\langle J\rangle$ is differentiable with respect to $s.$ 
In \emph{uniformly hyperbolic} systems, a stationary probability distribution, known as the SRB measure \cite{young}, exists, with respect to which ergodic averages converge starting from an open set in $\mathbb{R}^d$ containing the attractor. Under certain smoothness conditions on the map, the SRB measure is differentiable with respect to parameters \cite{ruelle2}, for small, smooth parameter perturbations. The assumption of uniform 
hyperbolicity is involved in the shadowing algorithm (section \ref{sec:nilss}) and in the data assimilation scheme (section \ref{sec:da}). There is a wealth of numerical and experimental evidence \cite{Gallavotti2006, angxiu2} that shows that physical systems exhibit \emph{quasi-hyperbolic} behavior. Hyperbolicity of the time-delayed system 
we consider in this paper, has been numerically verified by Huhn and Magri  \cite{Huhn2020_jfm} for a range of design parameters (numerical experiments are also 
presented in Figure \ref{fig:mean_angles} later in this paper).

\subsection{Tangent dynamics}
\label{sec:decomposition}
The tangent equation describes the response 
of the system's state to infinitesimal perturbations in a 
parameter $s \in \mathcal{S}$ in a neighborhood of a reference trajectory $\left\{u_n\right\}$. By introducing the shorthand $v_n := \partial_s u_n$,
the tangent equation is
\begin{align}
\label{eqn:tangent}
		v_{n+1} &= \partial_s f(u_n, \mathcal{S}) + (D_u f)(u_n, \mathcal{S})\; v_n \\
\notag v_0 &= 0 \in \mathbb{R}^d,
\end{align}
where $D_u$
denotes the differentiation with respect to the state vector.
We refer to the 
solutions $v_n$ as the inhomogeneous tangent solutions.
On setting the \emph{source} term in Eq. \ref{eqn:tangent} to zero, and starting with a 
non-zero initial perturbation, we obtain the time evolution of the perturbations in the initial state, denoted $q_n$. We
refer to $q_n$ as the homogeneous tangent solution whose time evolution is given by
\begin{align}
\label{eqn:homotangent}
		q_{n+1} = (D_u f)(u_n,\mathcal{S}) \; q_n,  
\end{align}
The solution $q_n$ is the derivative: $q_n := (D_u f_n)(u_0,\mathcal{S}) q_0$, which means that the homogeneous tangent equation is an iterative application of the chain rule. The homogeneous tangent solution is the difference between two orbits of $f$ at $n$, which are separated by an infinitesimal distance along $q_0$ at time $0.$
The inhomogeneous 
tangent solution, on the other hand, is the difference between two orbits of $f(\cdot, \mathcal{S})$ at infinitesimally different $s$, starting from the same initial condition. We can write 
down the following difference approximation of the 
inhomogeneous tangent equation:
\begin{align}
\label{eqn:finiteDifferenceInhomoTangent}
		v_n \approx \dfrac{f_n(u, s + \epsilon) - u_n}{\epsilon}.
\end{align}
The homogeneous tangent equation can also be approximately computed 
by finite differences
\begin{align}
    \label{eqn:finiteDifferenceHomoTangent}
    q_n \approx \dfrac{f_n(u + \epsilon q_0, \mathcal{S}) - u_n}{\epsilon},
\end{align}
Both these finite difference approximations are valid only up to an index $n$ for 
which the perturbed trajectory, and the original trajectory $\left\{u_n\right\}$, remain near each other. Because of chaos, for almost every direction $q_0$, the perturbed and unperturbed trajectories, exponentially diverge from each other. The 
finite difference approximations are 
bounded by $D/\epsilon,$ where $D$ is 
a scalar upper bound for the attractor within 
which all state vectors lie.
On the other hand, the tangent solutions $v_n$ and $q_n$, which are limits as 
$\epsilon \to 0$ of the right hand sides of Eq. \ref{eqn:finiteDifferenceInhomoTangent} and Eq. \ref{eqn:finiteDifferenceHomoTangent}, respectively, continue to diverge 
exponentially with $n$ in a chaotic system, unlike the finite difference approximations. 
That is, for large $n$, and almost every $q_0$,
$\norm{v_n}, \norm{q_n} \sim   e^{\lambda_1 n},$ where $\lambda_1 > 0$ is the largest characteristic \emph{Lyapunov exponent} \cite{Arnold1990}. 
In this paper, we refer to $1/\lambda_1$ as the \emph{Lyapunov time},  which is a timescale for the number of iterations needed to increase the norm of a linear perturbation by a factor $e.$ However, not all infinitesimal perturbations diverge exponentially. There are initial conditions, $q_0$, at every $u$ on a chaotic attractor that generate (asymptotically) exponentially decaying homogeneous tangent solutions along the orbit of $u$. This is because the space of tangent solutions $\mathbb{R}^d$, has the direct sum decomposition $\mathbb{R}^d  = E^u(u)\oplus E^s(u)\oplus E^c(u).$ The linear subspaces $E^u(u)$ and $E^s(u)$ contain tangent vectors at $u$ that exponentially grow/decay asymptotically under the tangent dynamics in Eq. \ref{eqn:tangent}, respectively. There is a maximum of $d$ possible Lyapunov exponents that give the asymptotic exponential growth/decay rates of tangent vectors (see \cite{Arnold1990} for Oseledets theorem). Hereafter, we assume that there are $d_u$ strictly positive Lyapunov exponents, which means that the unstable subspace at each point is $d_u$-dimensional. In a chaotic system, $d_u \geq 1$, i.e., the unstable subspace is at least one-dimensional at every $u.$ In this paper, the center subspace $E^c(u)$ consists of all the tangent vectors that asymptotically neither grow nor decay on an exponential scale, i.e. the tangent vectors in this subspace have a zero Lyapunov exponent. For example, consider the tangent vector $\mathcal{F}(u)$, whose flow is our primal system (Eq. \ref{eqn:ode}). 
If $f$ is a numerical discretization of the dynamical system in Eq. \ref{eqn:ode}, $\mathcal{F}$ approximately satisfies Eq. \ref{eqn:tangent} (it exactly satisfies the continuous-in-time formulation of Eq. \ref{eqn:tangent}). In this paper, we assume that $E^c(u)$ 
is one-dimensional at every $u$ and is spanned by a bounded tangent vector field, say $\tilde{\mathcal{F}}$, which exactly satisfies Eq. \ref{eqn:tangent}, 
\begin{align}
    \tilde{\mathcal{F}}(u_{n+1}) = (\partial_s f)(u_n,\mathcal{S}) + (D_u f)(u_n, \mathcal{S})\: \tilde{\mathcal{F}}(u_n),
\end{align}
The vector field 
$\tilde{\mathcal{F}}$ is approximated by the known vector field $\mathcal{F}$, and will be referred to as the center direction.
Note that this is a slight generalization of uniform hyperbolicity ($E^c$ is technically absent in a uniformly hyperbolic system), for which we assume the uniqueness and 
differentiability of the SRB measure \cite{pugh}.

\subsection{Adjoint dynamics}
\label{sec:adjointDynamics}
Exponential divergence also holds for adjoint equations starting from almost every initial condition. Fixing $N \in \mathbb{N}$, the 
$N$-time average $\langle J\rangle_N$ is affected by the primal solution 
at each $n\leq N$. If $u_n$ is infinitesimally perturbed, then, $J(u_m)$ is modified for all $m \geq n,$ causing the $N$-time average $\langle J\rangle_N$ to be altered. Viewed in this manner, at a fixed $\mathcal{S}$, $\langle J\rangle_N$ is a function of $N$ variables, $\left\{ u_n\right\}_{n=0}^{N-1},$ where each variable $u_m$ is, in turn, a function of $u_n,$ $n< m.$ That is, $\langle{J}\rangle_N \equiv 
\langle J\rangle_N(u_0, u_1, u_1, \cdots, u_{N-1}),$ with $u_n = f(u_{n-1},\mathcal{S})$. The adjoint solution at time $n$ is the response of $\langle J\rangle_N$ to an infinitesimal perturbation in $u_n,$ keeping the states prior to $n$ fixed at a reference orbit. At $n\leq N$, the adjoint solution is defined as 
\begin{align}
    \label{eqn:definitionadjoint}
    v^*_n := (D_{u_n} \langle J\rangle_N)^T(\left\{u_n\right\}_{n=0}^{N-1}) \in \mathbb{R}^d,
\end{align}
where $D_{u_n}$ refers to the total derivative with respect to $u_n$ and $^T$ stands for transpose. For $n =N$, Eq. \ref{eqn:definitionadjoint} gives
$v^*_N = (1/N) \: (D_u J)^T(u_N,\mathcal{S})$. By applying the chain rule, the adjoint vectors, $v^*_n,$ satisfy 
the inhomogenous 
adjoint equation, which is an iterative equation
\begin{align}
\notag
v^*_n &= (D_uf)^T(u_n,\mathcal{S})\; v^*_{n+1} + \dfrac{1}{N} (D_u J)^T(u_n)  \\
v^*_{N+1} &= 0.
\label{eqn:inhomoAdjoint}
\end{align}
The inhomogenous adjoint equation is solved backward in time 
with the zero vector as the initial condition at $N+1.$ 
The homogeneous adjoint solution is defined by setting the source term to zero, which yields
\begin{align}
q^*_n &= (D_uf)^T(u_n,\mathcal{S})\; q^*_{n+1}.
\label{eqn:homoAdjoint}
\end{align}
This equation, which is solved backward in time with a non-zero 
initial condition at $N,$ is an 
iterative application 
of the chain rule, at a fixed $\mathcal{S}$,
to the definition 
$q^*_n := (D_u (f_{N-n}\cdot q_N^*))^T(u_n)$. 
In other words, the homogeneous adjoint solution at time $n$ is the sensitivity to $u_n$ of 
the solution at time $N$ projected along $q^*_N$, $N \geq n.$
Similar to the tangent solutions, 
the adjoint solutions asymptotically 
grow exponentially, backward in time, 
for almost every initial condition, i.e., for large 
$N$, $q^*_0, v^*_0 \sim {\cal O}(e^{\lambda_1 N}).$ Intuitively, we can understand this growth as complementary to the growth of tangent solutions. That is, since infinitesimal perturbations to the state increase in norm exponentially forward in time, we expect the solution at a given time to be more sensitive (exponentially) to its far past compared to its recent past. The Lyapunov exponents 
characterizing the adjoint dynamics' growth are the same as those for the tangent equations. The tangent and adjoint solutions are connected by bi-orthogonality, which means that, homogeneous adjoint solution and 
tangent solutions that are associated with two different Lyapunov exponents, are orthogonal to each other. In particular, the tangent solutions with negative LEs, which span $E^s$, must be orthogonal to adjoint solutions associated to positive or zero LEs. In other words, it can be shown that the unstable adjoint subspace, consisting of all the adjoint solutions that exponentially 
grow (at most at $d_u$ different asymptotic rates) backward in time, is orthogonal to 
$E^s$ and $E^c.$ Similarly, the stable adjoint subspace is orthogonal to $E^u$ and $E^c$ \cite{Kuptsov2012}. 

\subsection{Automatic differentiation to compute tangent and adjoint solutions}
\label{sec:adintro}
We consider another class of linear perturbation 
methods: automatic differentiation (AD). 
Given a program, with \emph{output} $O$ and \emph{input} $I$,  where $O$ and $I$ can be 
scalar or vector-valued, AD obtains the derivative 
$d O/d I$. In forward mode, the program is traversed 
sequentially, and each line is differentiated with respect to $I$ using 
the derivatives of the variables computed in the previous lines. Ultimately, 
$dO/dI$ is obtained exactly. In reverse mode AD, the derivative $d O/d I$ is again 
obtained exactly, but by traversing the program in reverse order and using the chain rule to update the derivatives. 

As noted in section \ref{sec:decomposition} and section \ref{sec:adjointDynamics}, each of the four linear perturbation solutions 
discussed (homogeneous/inhomogeneous tangent/adjoint solutions) can be written in the form of a derivative. By specifying $O$ and $I$ appropriately, 
all of them can be computed through AD.
 From their derivative-based definitions, the functions that specify $O$ only require the primal solver (i.e., a time-integrator, $f(\cdot, \mathcal{S})$) and the definition of the objective function. As a result, in order to compute linear perturbation solutions, AD does not require the user to compute the Jacobian, and eliminates the need for tangent/adjoint calculations through their respective iterative equations. For example, consider the AD solution of the inhomogeneous adjoint equation (Eq. \ref{eqn:inhomoAdjoint}). To compute it, we define a function that returns a value $O = (1/N)\sum_{n=0}^{N-1} J(f_n(u_0), \mathcal{S})$. The function takes the argument $I = u_0$ time-integrates the primal for $N$ steps, evaluates $J$ at each step for averaging at the end. Then, $dO/dI$ computed by applying AD in reverse-mode is the solution of the inhomogeneous adjoint equation at time 0. The reader is referred to texts on AD (e.g. Ch 3 \cite{griewank3}, Ch. 10 \cite{griewank10} and Ch. 15 \cite{griewank15} of Griewank and Walther) for the application of AD to replace tangent/adjoint solvers. In this paper, our focus is the implementation of the shadowing algorithm \cite{angxiu1} to compute sensitivities. We use AD to replace the tangent/adjoint solvers needed within the shadowing algorithm. The inputs and outputs to AD must be defined appropriately in the AD version of shadowing, as discussed in \ref{sec:appxAD}.
 
While development time is reduced by AD through the elimination of hand-differentiation, 
AD has compile-time and run-time overheads, which depend on the AD software used.
In this paper, we use the AD package \verb+Zygote.jl+ \cite{Innes} in Julia, which uses the language's multiple dispatch feature to compute derivatives. Several AD library options exist in languages popular in scientific computing \cite{autodiff}. Some combine modern language-level features (e.g. multiple dispatch or operator overloading on dynamic types in Julia) with algorithmic advances \cite{InnesEdelman} to achieve time and memory efficiency \cite{ChrisComparison} when compared to traditional solvers for perturbation equations, (see e.g. optimization in PerforAD in Python \cite{PerforAD}).

\begin{figure}
    \centering
    \includegraphics[width=0.8\textwidth]{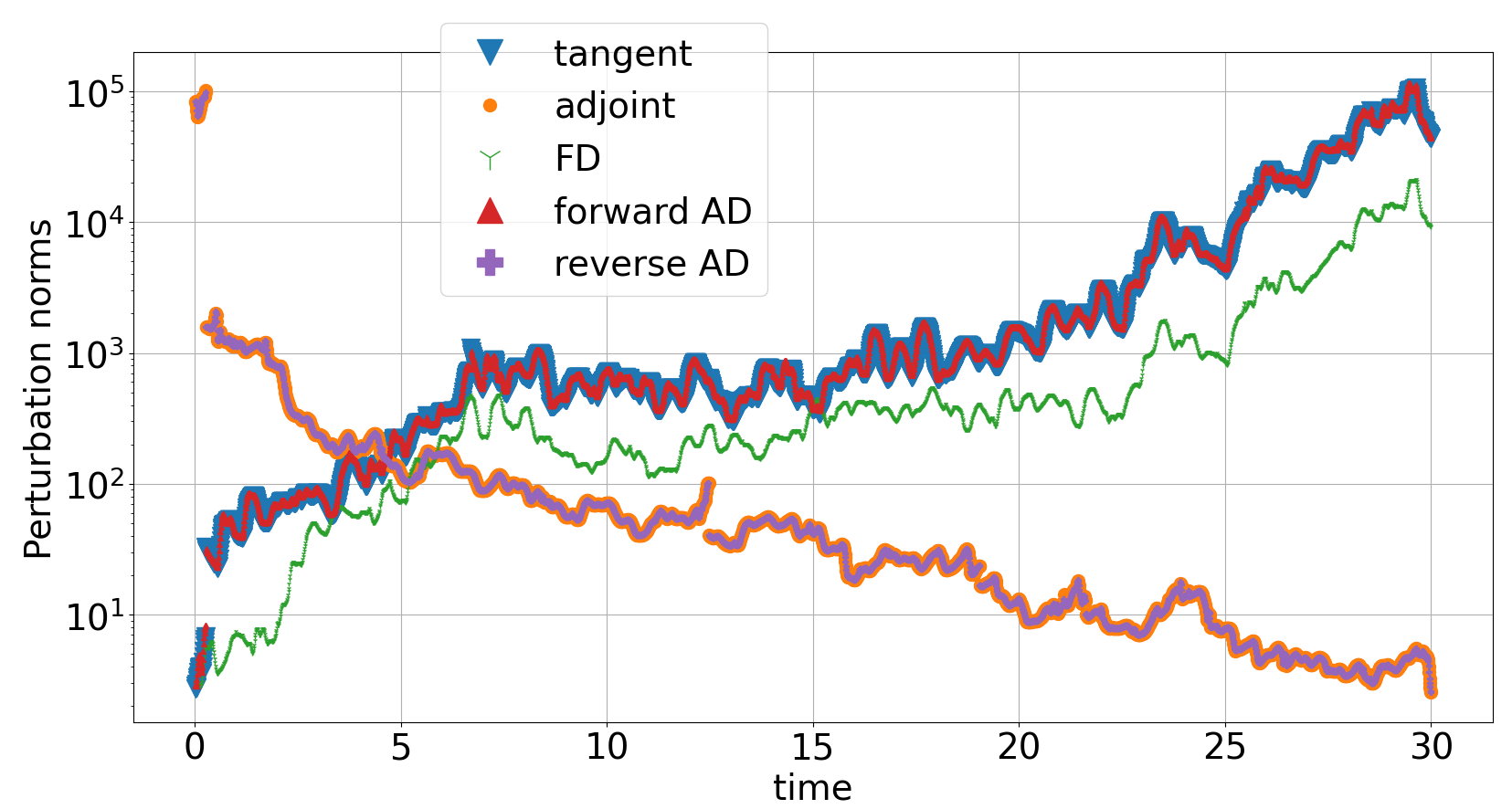}
    \caption{$l^2$ norms of the perturbation vectors computed through the homogeneous tangent (blue triangle), adjoint equations (orange circle), finite difference (green Y), forward-mode AD (red triangle) and reverse-mode AD (purple plus), are shown as a function of time, for the time-delayed model of section \ref{sec:rijke}.}
    \label{fig:perts}
\end{figure}

In Figure \ref{fig:perts}, we plot the $l^2$ norms of linear perturbations computed using the four different methods discussed. The primal system that supplies $f$ is the chaotic acoustic model (section \ref{sec:rijke}). All of the linear perturbation methods evolve with exponentially increasing norms. As noted earlier, finite difference saturates on the order ${\cal O}(1/\epsilon)$ since the attractor is bounded and, therefore, so is the norm of the difference between any two solutions. The finite difference results shown as green Y's in Figure \ref{fig:perts} are calculated with an initial perturbation with norm $10^{-4}.$ The norm of the finite difference increases exponentially before saturating at about $10^4$. Forward-mode AD results, which compute the tangent solutions exactly, closely approximate the latter. The tangent solutions computed using Eq. \ref{eqn:homotangent} and forward-mode AD, and the adjoint solution computed using Eq. \ref{eqn:homoAdjoint} and reverse-mode AD, show unbounded exponential growth. The slopes of the perturbations on the logarithmic scale ($\approx 0.2$) indicate the largest Lyapunov exponent of the chaotic acoustic model (section \ref{sec:rijke}). This is a manifestation of the \emph{butterfly effect}. Next, we explain how this effect leads to the breakdown of traditional sensitivity algorithms in chaotic problems.  
\subsection{The problem with computing sensitivities of ergodic averages using conventional methods}
To compute the sensitivity of a time-averaged quantity, 
one could potentially use the tangent, adjoint solutions, 
or finite difference, or AD. Consider the problem of sensitivity computation of a finite-time 
average $\langle J\rangle_N.$ Using the tangent or adjoint solutions solved iteratively using Eq. \ref{eqn:tangent} and Eq. \ref{eqn:inhomoAdjoint} (or computed using AD) respectively, one can obtain the derivative of $\langle J\rangle_N$ as
\begin{align}
    \label{eqn:tangentSensitivity}
    d_s \langle J\rangle_N (u_0, s) &= 
    \dfrac{1}{N}\sum_{n=0}^{N-1} x^*_n\cdot v_n \\
    \label{eqn:adjointSensitivity}
    &= \dfrac{1}{N}\sum_{n=0}^{N-1} x_n\cdot v^*_n,
\end{align}
where, for notational convenience, we have defined  
$x_n := (\partial_s f)(u_{n-1}, \mathcal{S}),$ and $x_n^* := (D_u J)^T(u_n, \mathcal{S}).$
Both equations, which can be derived using the chain rule, along with their AD counterparts discussed in section \ref{sec:adintro}, are standard 
in sensitivity analysis. To reap computational benefits, Eq. \ref{eqn:tangentSensitivity}, or forward-mode AD, is used when the set of 
observables has a larger dimension than the parameter space. The tangent 
solution $\left\{v_n\right\}$ corresponding to a parameter, can be used to compute the 
derivatives of time averages of 
all the observables with respect to that parameter. By contrast, 
when the number of parameters exceeds the number of observables, Eq. \ref{eqn:adjointSensitivity}, or reverse-mode AD is the preferred approach to compute sensitivities 
since the same sequence $\left\{ v^*_n\right\}$ computed for a given $\langle J\rangle_N$ is used to compute the 
gradient of $\langle J\rangle_N$ with respect to all the parameters in $\mathcal{S}.$

In chaotic systems, the above approach yields values of $d_s \langle J\rangle_N$ 
that exponentially increase with $N$, as $N\to\infty.$ 
However, the quantity of 
interest, $d_s \langle J\rangle$ (Eq. \ref{eqn:qoi}), in which the limit $N\to \infty$ is taken before the derivative with respect to $s$, is bounded. Hence the derivative of the ergodic average $d_s \langle J\rangle$ is not the same as the derivative of the finite-time average $d_s \langle J\rangle_N$ in the limit $N\to \infty.$ Thus, conventional methods 
for sensitivity computation are not applicable to the computation of derivatives 
of ergodic averages in chaotic systems. One early approach to circumvent this problem is the ensemble sensitivity method \cite{Lea2000, Eyink2004} in which $d_s \langle J\rangle$ is approximated by a sample average of sensitivities 
computed by using Eq. \ref{eqn:tangentSensitivity} or Eq. \ref{eqn:adjointSensitivity} over a small $N$ (comparable to one Lyapunov time). The accuracy of this method improves as $N$ increases, provided that the number of samples increases exponentially with $N$. This makes the method prohibitively expensive in practice \cite{Eyink2004, nisha-ens}. In the next section, we describe Non-Intrusive Least Squares Shadowing (NILSS) due to Ni \emph{et al.} \cite{angxiu1}, which is a more efficient approach for computing the same quantity.

\subsection{Non-Intrusive Least Squares Shadowing}
\label{sec:nilssidea}
Owing to the \emph{shadowing lemma} \cite{shadowing}
for uniformly hyperbolic systems, it has been shown by Wang \cite{qiqi} that
there exists a unique perturbation direction $v^{\rm sh}$ -- the tangent shadowing perturbation -- for which
the tangent equation (Eq. \ref{eqn:tangent}) has a bounded solution for all time.
The shadowing perturbation $v^{\rm sh}$ is an inhomogeneous tangent solution (i.e., $v^{\rm sh}$ satisfies Eq. \ref{eqn:tangent}). However, unlike the conventional tangent solution, $v$, it does not exhibit an unstable growth. This constraint is 
used by the NILSS algorithm \cite{angxiu1} to solve for $v^{\rm sh}$ over a long, but finite-time, duration. In particular, NILSS \cite{angxiu1} constructs an approximation of the shadowing perturbation by subtracting from $v$ an unstable tangent vector at every point along a trajectory. The unstable tangent vector to be subtracted is represented in an orthonormal basis of the unstable tangent subspace ($E^u$). The orthonormal basis, is in turn, computed by propagating at least 
as many tangent vectors as the dimension of the unstable subspace, under the 
homogeneous tangent dynamics (Eq. \ref{eqn:homotangent}), along with repeated 
normalization. This procedure is typically used in the computation of Lyapunov vectors \cite{Kuptsov2012, Ginelli2007}. 
\\

Let $Q_n$ be a $d\times d_u$ matrix whose columns form an 
orthonormal basis of $E^u(u_n)$. In the NILSS algorithm \cite{angxiu1}, 
the total time duration $N$ is divided into multiple short time segments, \emph{checkpoints}, such that each short segment is comparable to the Lyapunov time. We shall simplify the setting by considering a time segment to be one timestep (i.e., every timestep is a checkpoint); we delay a discussion on this simplification until the end of this section. The tangent shadowing perturbation can be expressed as
\begin{align}
v^{\rm sh}_n = v_n + Q_n \: a_n, \; 1\leq n \leq N,
\label{eqn:tangentShadowing}
\end{align}
where $a_n \in \mathbb{R}^{d_u}$ is a vector of coefficients. NILSS computes the sequence of vectors $\left\{a_n\right\}$ along the trajectory by solving a minimization problem for the norms of $\left\{v^{\rm sh}_n\right\}.$ For a complete description of the NILSS algorithm, the reader is referred to \cite{angxiu1}, where the derivation 
of the algorithm for the time-continuous case is presented. Since perturbations along the center direction neither grow nor decay exponentially, in NILSS, the center direction is  excluded from $Q$. Its effect on the sensitivity is added back later (see section 2 of \cite{angxiu1}). When the map $f$ is a  time-discretized ODE, it has a center direction corresponding (but not exactly equal) to the center direction of the ODE. Thus, we also take into account, in the discrete algorithm, modifications due to this center direction, when $f$ models ODEs, as we shall see in section \ref{sec:center}.

The NILSS problem minimizes the norms of the shadowing perturbation sequence $\left\{ v^{\rm sh}_n\right\}$. The Lagrangian of this optimization problem can be written as 
\begin{align}
		{\cal L}^{\rm sh}(\left\{a_n\right\}, 
		\left\{\beta_n\right\}) &:= \sum_{n=1}^N \|v^{\rm sh}_n\|^2 + 
		\sum_{n=1}^{N-1} \beta_n (a_{n+1} - R_{n+1} a_n - \pi_n) \\
		\label{eqn:nilss_lsp}
		&= \sum_{n=1}^N \Big( \|v_n\|^2 + \|a_n\|^2 + 2 v_n \cdot Q_n a_n\Big) +  \sum_{n=1}^{N-1} \beta_n (a_{n+1} - R_{n+1} a_n - \pi_n),
\end{align}
where Eq. \ref{eqn:nilss_lsp} uses Eq. \ref{eqn:tangentShadowing}, and the fact that 
$Q_n^T Q_n$ is the $d_u\times d_u$ identity matrix. Here, $\left\{\beta_n\right\}$ is a sequence of Lagrange multipliers that imposes a sequence of equality constraints at every timestep to ensure the continuity of the shadowing perturbation (section \ref{sec:lss}). We solve the above problem to obtain a sequence $\left\{a_n\right\}$, and then, to obtain the shadowing perturbation through Eq. \ref{eqn:tangentShadowing}. Subsequently, we compute the required sensitivity through Eq. \ref{eqn:tangentSensitivity}, 
with the (exponentially growing) tangent solution, $v_n$, replaced with the 
shadowing tangent solution $v^{\rm sh}_n.$ This yields
\begin{align}
\label{eqn:tangentNILSSsensitivity}
    d_s \langle J\rangle \approx \frac{1}{N}\sum_{n=0}^{N-1} x^*_n \cdot 
    v^{\rm sh}_n,
\end{align}
as shown in Appendix C of \cite{angxiu1}, or in 
Theorem LSS of \cite{qiqi}. The 
same shadowing perturbation $v^{\rm sh}_n$ is used to compute the sensitivities with respect to all $J \in \mathcal{J}.$ On the other hand, the tangent NILSS algorithm has to be repeated for every parameter $s \in \mathcal{S}$ in order to compute the corresponding shadowing perturbations. When the parameter space is higher-dimensional when compared to the observable space, the adjoint version of NILSS is preferred.

\subsubsection{Adjoint non-intrusive least squares shadowing}
While the adjoint algorithm can be obtained via reverse-mode automatic differentiation of
tangent NILSS, the theoretical basis for the existence of an adjoint shadowing 
perturbation is developed in \cite{angxiu-siam}. The adjoint algorithm, known as Non-Intrusive Least Squares Adjoint Shadowing (NILSAS), is presented in \cite{angxiu-nilsas} with an application to a fluid flow problem in \cite{angxiu2}. Here, we focus on the discrete time case, and for simplicity, each time segment corresponds to one iteration
of the map $f.$ Analogous to tangent NILSS, in the adjoint version, an 
adjoint shadowing perturbation $v^{\rm sh *}_n$ is computed, which solves the 
inhomogeneous adjoint equation (Eq. \ref{eqn:inhomoAdjoint}). Extending the analogy further, an unstable adjoint vector is subtracted from the conventional adjoint solution to obtain $v^{\rm sh *}_n$ as 
\begin{align}
    \label{eqn:adjointShadowing}
    v^{\rm sh *}_n = v^*_n + Q^*_n a^*_n,
\end{align}
where $Q^*_n$ is an orthonormal basis for the unstable adjoint subspace, $(E^s \oplus E^c)^{\perp}$, and $a^*_n \in \mathbb{R}^{d_u}$ is a set of coefficients. The orthonormal basis $Q^*_n \in \mathbb{R}^{d\times d_u}$
is achieved by 
iterating at least $d_u$ adjoint vectors backward in time, using 
Eq. \ref{eqn:homoAdjoint}, repeatedly normalizing with QR factorization. 
That is, the orthonormalization procedure is identical to that for $\left\{ Q_n\right\}$, but using the sequence of Jacobian transposes,
instead of the Jacobians, and with time-reversal. The particular set of coefficients $a^*_n$ needed to find the adjoint 
shadowing sequence $v^{\rm sh *}_n$, which is a bounded solution of the inhomogeneous adjoint equation, is found as a solution of a least squares problem, which is also analogous to Eq. \ref{eqn:nilss_lsp}. In order to compute the quantity of interest, we replace the conventional (exponentially growing) 
inhomogeneous adjoint solution in Eq. \ref{eqn:adjointSensitivity} with 
the adjoint shadowing solution
\begin{align}
    \label{eqn:adjointNILSSsensitivity}
    d_s \langle J\rangle \approx \dfrac{1}{N}\sum_{n=0}^{N-1} v^{\rm sh *}_n \cdot 
    x_n.
\end{align}
We now comment on the error-vs-cost trade-off of the tangent and 
adjoint shadowing algorithms, noting that a careful analysis of this trade-off is problem-specific and beyond the scope of this paper. In both algorithms, the most expensive computation in the shadowing algorithm is the solution of the minimization problem. The size of the minimization 
problem is directly proportional to $d_u$ and $N.$ For the same overall time duration 
$N,$ choosing a larger time segment between checkpoints reduces the size of the minimization problem, since the QR factorization and the equality constraints that enter the problem (Eq. \ref{eqn:nilss_lsp}) are only computed at the boundaries of the time segments. 

For example, in the time delay acoustic problem of this paper, we could theoretically choose a segment size over which the linear perturbations increase in norm by, say, a factor of 2. From Figure \ref{fig:perts}, such a segment size is about 1 time unit (100 timesteps). Then, we need to 
perform QR factorization only every 100 timesteps. While the QR factorization in itself is not the computational bottleneck, the size of the least squares problem shrinks by a factor of 100, as compared to checkpointing every timestep. In spite of the additional cost, we choose to checkpoint every timestep in the considered problem. However, we choose a size $N$ such that the ${\cal O}((N d_u)^3)$ calculation of the least squares solution is neither memory-constrained nor is a prohibitive computational expense. Then, we repeat the shadowing algorithm $M$ times and 
sample-average the shadowing sensitivities. This procedure is effectively the same as computing the sensitivity of an $MN$-time average by executing the shadowing algorithm once, provided that $N$ is large enough for the convergence of the ergodic averages. The reason for sample-averaging shadowing sensitivities as opposed to segmenting a long shadowing algorithm is that we observe, for the time-delayed acoustic model, higher condition numbers of the least squares problem. Thus, we take the approach of checkpointing every timestep, and solving smaller least squares problems. The size of each problem is chosen large enough for ergodic averages to converge while also curtailing the computed shadowing perturbations (from the solution of the least squares problem) to an ${\cal O}(1)$ norm, 
at every timestep.  The advantage of this approach is the simpler program whose computational cost is nearly the same as the $MN$-sized checkpointed NILSS, but produces better-behaved shadowing perturbations. 

\subsection{Error in shadowing and alternatives}
The NILSS algorithm is not guaranteed to converge to the true value of the 
sensitivity, as $N \to \infty$. To see why, we first note that the 
shadowing sensitivities computed by NILSS, which are the right hand sides of 
Eq. \ref{eqn:tangentNILSSsensitivity} and Eq. \ref{eqn:adjointNILSSsensitivity}, 
are ergodic averages along a true orbit of a system at an infinitesimally perturbed $s$. For a mathematically rigorous explanation, 
the reader is referred to \cite{qiqi}; a qualitative explanation is also 
included in section \ref{sec:da}. Now, due to a perturbation 
in $s,$ the stationary probability distribution on the attractor is perturbed as well, but this 
perturbation is excluded by NILSS. In general, the
ergodic averages along a shadowing orbit (a true orbit) of a perturbed system 
do not converge to the expectation with respect to the 
stationary distribution of the unperturbed system. Hence, there is a systematic 
error in NILSS, which has been recently studied in \cite{angxiu-shadowing}, 
along with the corrections that can be made to reduce the error.
In view of this shadowing error, we must mention that a few alternatives have recently appeared. In particular, the space-split sensitivity 
method \cite{nisha-s3} is an ergodic-averaging method to compute 
Ruelle's linear response formula \cite{ruelle2}, which specifies the required sensitivity exactly.
However, based on the preliminary formulation \cite{nisha-s3, adam}, 
it is more complex to implement than the shadowing method. Thus, in 
applications where a systematic error in the computed sensitivities is 
not a serious impairment, such as in the parameter estimation and data assimilation  problems considered in sections \ref{sec:optim} and \ref{sec:da}, respectively, shadowing methods may be preferred. Another sensitivity computation method, also based on    
Ruelle's linear response formula \cite{ruelle2}, is known as blended response \cite{Abramov_2007}, in which short-term and long-term responses are computed 
using different methods. For the long-term response to unstable perturbations, 
Ruelle's formula, which is exact, is approximated using a Fluctuation-Dissipation 
theorem for non-equilibrium settings \cite{lucarini}. However, this approach is computationally expensive and also inexact. Thus, we focus on the 
shadowing-based methods for computing sensitivities in this paper. Moreover, shadowing methods have successfully been applied to dissipative models in fluid mechanics \cite{angxiu2, Huhn2020_jfm}, as relevant to this study.

\section{Tangent/adjoint shadowing algorithm}
\label{sec:nilss}
We provide a step-by-step description of the tangent and adjoint NILSS algorithms. The reader is referred 
to \cite{angxiu1} and \cite{angxiu-nilsas} for the original descriptions of tangent and adjoint NILSS, respectively. We consider the discrete algorithm without the checkpointing scheme. We also adopt a simplified presentation for which 
the same program can be used for implementing both tangent and adjoint NILSS, with a minimal modification. Further, here, we introduce AD to compute the needed 
tangent and adjoint solutions. This automatic-differentiated unified program for the tangent/adjoint NILSS, shall be referred to as the 
AD shadowing algorithm. 

\subsection{Shadowing algorithm: inputs and outputs}
We present a shadowing algorithm that 
constructs a sequence of shadowing perturbations 
$v^{\rm sh}_n, 1\leq n\leq N.$ The algorithm 
takes as inputs, a sequence of $d\times d$ matrices $A_n$, 
and a sequence of $d$-length vectors, $b_n,$ to return 
tangent or adjoint shadowing sensitivities (defined in 
Eq. \ref{eqn:tangentNILSSsensitivity} and Eq. \ref{eqn:adjointNILSSsensitivity}, respectively)
in the following two scenarios. 
\begin{itemize}
    \item \textbf{Case 1 (tangent):} The sequence $\left\{ A_n\right\}$ is set to the Jacobian 
    matrix sequence $\left\{ D_u f(u_n, \mathcal{S})\right\}$ along a reference trajectory $\left\{u_n\right\}$. The sequence $b_n$ is the parameter perturbation at $n$, 
i.e., $b_n = x_n.$ Then, the algorithm returns 
the sequence of tangent shadowing perturbation vectors
at $u_1, u_2,\cdots, u_N$, namely, $\vsh_1, \vsh_2, 
\cdots , \vsh_N.$ These shadowing perturbations can be used to approximately compute the $l$ sensitivities, 
$d\langle \mathcal{J}\rangle/ds.$
\item \textbf{Case 2 (adjoint):} Now, on the other hand, 
defining $n' := N + 1 - n$ 
suppose $A_n$ is the transpose of the
Jacobian matrix at $n'$, that is,
$A_n = (D_u f)^T(u_{n'})$. Set 
$b_n := x^*_{n'+1}$. 
In this case, the 
algorithm returns the sequence of adjoint shadowing  
perturbation vectors at $u_{N+1}, u_N, u_{N-1}, \cdots, u_2$, 
namely, $\vsh_1, \vsh_2, \cdots, \vsh_N$ (in the unified presentation, we drop the superscript ``*'' used for adjoint solutions). 
The adjoint shadowing perturbation sequence is used to 
compute the $p$ sensitivities, $D_{\mathcal S} \langle J\rangle.$
\end{itemize}
 A time reversal is accomplished in adjoint shadowing simply by reversing the indexing of the 
  input sequences; the output sequence of adjoint shadowing perturbations is obtained in time-reversed order. We 
use the term \emph{shadowing perturbation} to refer
to both tangent and adjoint shadowing perturbations. Note that $u_0$ must be a point on the attractor sampled according to the stationary probability distribution on the attractor. That is, the primal system must be simulated for a run-up time long enough for time-averages to converge. A solution state obtained after such a run-up time is chosen as $u_0$.
\subsection{Evolution of homogeneous and inhomogeneous 
perturbations with repeated normalization}
\label{sec:nloop}
Our goal is to compute the tangent or adjoint 
shadowing perturbation using Eqs. 
\ref{eqn:tangentShadowing} and \ref{eqn:adjointShadowing}, 
respectively. Toward this goal,
we solve i) at least $d_u$ homogeneous 
tangent equations, or $d_u$ homogeneous adjoint equations, and ii) $p$ inhomogeneous 
tangent equations or $l$ inhomogeneous adjoint equations. 
The common form of the homogeneous equation, which amounts 
to solving the tangent equation forward (in case 1) 
or the adjoint equation backward 
in time (in case 2), is given by
\begin{align}
    \label{eqn:homogeneousPerturbation}
    q_n^i = A_{n-1} q_{n-1}^i, \;\; n = 1,\cdots, N, \;\; 
    1\leq i\leq d_u.
\end{align}
We define $Q_n$ to be an $n\times d_u$ matrix with columns $q_n^i.$ The following equation gives the evolution of the inhomogeneous 
tangent solution forward in time in case 1, and the inhomogeneous adjoint 
solution backward in time in case 2,
\begin{align}
    \label{eqn:inhomgeneousPerturbation}
    v_n = A_{n-1} v_{n-1} + b_n, \;\; n = 1,\cdots, N.
\end{align}
At each $n,$ we normalize both the homogeneous 
and inhomogeneous perturbations by QR factorization. We 
choose $d_u$ pseudo-random vectors in $\mathbb{R}^d$ as initial conditions  
$q_0^i$, $1\leq i\leq d_u$. The initial condition 
for Eq. \ref{eqn:inhomgeneousPerturbation}, $v_0$, is set to $0 \in \mathbb{R}^d.$ Beginning with $n=1$, we perform the following 
loop until $n = N$.
\begin{enumerate}
    \item Advance 
    Eq. \ref{eqn:homogeneousPerturbation} by one timestep 
    for each $1\leq i \leq d_u.$ This 
    can be written as $Q_n \longleftarrow A_{n-1} Q_{n-1}.$
    
    \item QR-factorize the matrix $Q_n$ and set $Q_{n}$ to the 
    obtained ``$Q$''. Let the ``$R$'' from QR factorization be stored 
    as $R_{n}$. Thus, each $q_n^i$, $1\leq i \leq d_u$ is now a 
    unit vector.
    
    \item Obtain $v_n$ from $v_{n-1}$ by advancing Eq. \ref{eqn:inhomgeneousPerturbation} by one timestep.
    
    \item Set $\pi_n := Q_n^T v_n$, which is a $d_u$-length vector of 
    orthogonal projections of $v_n$ along $q_n^i.$
    
    \item Project out the unstable components of $v_n$. That is, set $v_n \longleftarrow v_n - \pi_n Q_n.$
    
    \item Go to step 1 with $n\longleftarrow n + 1$, or stop if $n = N.$
 \end{enumerate}
First, using the above orthonormalization procedure, 
$Q_n$ converges to an orthonormal basis for 
the unstable tangent (adjoint) subspace in case 1 (case 2).
Secondly, in case 1, we note that the above $n$-loop must be executed only once for 
the sensitivity with respect to $s$ of all $J \in \mathcal{J}$. 
Similarly, in case 2, the $n$-loop must be called just once if 
we wish to compute the sensitivity of $\langle J\rangle$ with respect to all the parameters $\mathcal{S}$. In other words, to obtain $D_{\mathcal S} \mathcal{J}$, in case 1 (tangent shadowing), 
the $n$-loop must be run as many times as the number of parameters (= $p$), 
and in case 2 (adjoint shadowing),
as many times as the number of objective functions (= $l$). Thirdly, the sequence of matrices $R_n$ can be used to obtain 
the Lyapunov exponents. In particular, if the $k$th diagonal 
element of the matrix $R_n$ is written as $R_n^k$, then,
the $k$th Lyapunov exponent $\lambda_k \approx (1/N)
\sum_{n=0}^{N-1} \log{|R_n^k|}$. This can be easily seen  by recasting the definition of Lyapunov exponents (\cite{Arnold1990}) and as an ergodic average. 
\subsection{Minimizing the growth of the shadowing perturbation sequence}
\label{sec:lss}
At the end of the $n$-loop described in section \ref{sec:nloop},
we have at our disposal the following sequences of vectors or matrices, 
where at each $n \leq N+1$,
\begin{itemize}
    \item $v_n$ is the inhomogeneous perturbation orthogonalized with respect to the 
unstable tangent (adjoint) subspace in case 1 (case 2).
\item $\pi_n$ consists of the orthogonal projections (before 
the orthogonalization) of $v_n$ on the unstable tangent (adjoint) 
subspace in case 1 (case 2).
\item $Q_n$ is a $d\times d_u$ matrix that forms an orthonormal basis for the unstable 
tangent (adjoint) subspace at each $n$ ($n'$) in case 1 (case 2), and, 
\item $R_n$ is a $d_u \times d_u$ matrix that contains the one-step growth factors 
of $Q_n$ under the tangent (adjoint) dynamics in case 1 (case 2).
\end{itemize}
In practice, a finite spin-up time, typically on the order of Lyapunov time, is needed for the convergence of $Q_n$ to an orthonormal basis for the true unstable (tangent/adjoint) subspace. 
We can write the ansatz for the shadowing perturbation sequence (Eq. \ref{eqn:tangentShadowing} and Eq. \ref{eqn:adjointShadowing}) 
in a form that is applicable to both tangent and adjoint shadowing 
sequences, denoted here as $v^{\rm sh},$
\begin{align}
    \label{eqn:commonShadowing}
    v^{\rm sh}_n = v_n + Q_n a_n. 
\end{align}
Here the sequence $a_n$ is the unknown $d_u$-length vector, which we shall solve for. In case 2, 
the sequence $a_n$, and subsequently 
$v^{\rm sh}_n$, are obtained in time-reversed order by virtue of 
time-reversing the inputs $A_n, b_n$ to the $n-$loop. 
In particular, $v^{\rm sh}_n$ is the adjoint shadowing perturbation 
at time $N+2-n.$ In order to 
solve for $a_n,$ we start by multiplying 
Eq. \ref{eqn:commonShadowing} by $A_n$ and adding $b_{n+1}$ to both 
sides of the equation, 
\begin{align}
\label{eqn:anDerivationFirstStep}
    A_n v^{\rm sh}_n + b_{n+1} = A_n v_n + b_{n+1} + A_n Q_n a_n.
\end{align}
Since the shadowing perturbation solves Eq. \ref{eqn:inhomgeneousPerturbation}, the left hand side is $v^{\rm sh}_{n+1}.$ Using steps 3 to 5 of the $n$-loop in section \ref{sec:nloop}, 
the first two terms on the right hand side of Eq. \ref{eqn:anDerivationFirstStep} become $v_{n+1} 
+ Q_{n+1} \pi_{n+1}.$ From step 2 of the $n$-loop, $A_n Q_n = Q_{n+1} R_{n+1}.$ Thus, 
\begin{align}
\label{eqn:anRelation}
    v^{\rm sh}_{n+1} = v_{n+1} + Q_{n+1} \pi_{n+1} + Q_{n+1} R_{n+1} a_n.
\end{align}
From Eq. \ref{eqn:commonShadowing}, 
the left hand side of the above equation is also equal 
to $v_{n+1} + Q_{n+1}a_{n+1}.$ We obtain 
the following iterative relationship for $a_n,$ 
after multiplifying both sides by $Q_{n+1}^T$
\begin{align}
    \label{eqn:iterativea}
    a_{n+1} = \pi_{n+1} + R_{n+1} a_n.
\end{align}
This is the equality constraint that must be added to the NILSS problem, whose Lagrangian is in Eq. \ref{eqn:nilss_lsp}.
Hence, Eq. \ref{eqn:iterativea} is also one of the KKT conditions ($D_{\beta_n} \mathcal{L}^{\rm sh} = 0$) of the NILSS optimization problem. 
Although one can theoretically solve Eq. \ref{eqn:iterativea} starting from a random guess for $a_1$ and iterating, this does not provide accurate results in practice; Eq. \ref{eqn:iterativea} is not a well-conditioned problem for $\left\{a_n\right\}$ due to the exponential growth of the round-off errors in $\left\{R_n\right\}$, which tend to accumulate upon iteration. Thus, following \cite{angxiu1, angxiu-nilsas}, we resort 
to the direct method of solving for the entire sequence $\left\{ a_n\right\}$ at once (Appendix A of \cite{angxiu-nilsas} and \cite{angxiu1}). The direct method is to solve the following system of linear equations for $\left\{a_n\right\}$
\begin{align}
    G X = H, 
\end{align}
where 
\begin{itemize}
	\item 
$G$ is an $N d_u \times (N +1) d_u$ block matrix with 
$d_u\times d_u$ blocks given by
				\begin{align}
					G := \begin{bmatrix}
						-R_1 & I  & 0 & \cdots & 0 & 0   \\
						0 & -R_2 & I & \cdots & 0 & 0 \\
						0 & \cdots & \cdots & \cdots & 0 & 0 \\
				\cdots & \cdots & \cdots & \cdots & \cdots & \cdots \\
				0 & \cdots & \cdots & -R_{N-1} & I & 0 \\
							0 & \cdots & \cdots & \cdots & -R_{N} & I 
					\end{bmatrix},
				\end{align}
				where $I$ is the $d_u\times d_u$ identity matrix,
	\item 
$X$ is an $N\times d_u$ vector consisting of $[a_0,\cdots,a_{N}]$ and,
\item $H$ is an $N d_u$-length vector containing 
		the sequence $[\pi_1,\cdots,\pi_{N}].$ 
\end{itemize}
The solution of the underdetermined system that minimizes the norm of 
$X$ is given by 
$$ X = G^T (G G^T)^{-1} H.$$
\subsection{Modifications due to the center direction}
\label{sec:center}
Whether in tangent or adjoint shadowing, a better accuracy is obtained 
if the center direction, which is approximately 
the right hand of an ODE when $f$ is a time-discretization of the ODE, is given a special treatment. In tangent shadowing, $d_u$ can be set to the 
number of positive LEs plus 1, so that the center direction 
is treated as an unstable direction. However, in some problems, this may increase
the condition number of the least squares problem for $X$. This  
leads to a poorer minimization of $X,$ which in turn increases 
the norm of the shadowing perturbation, when compared to the following 
alternative. As suggested in \cite{angxiu1}, we project out the center 
components of both the homogeneous and inhomogeneous tangents 
and add the contribution to the sensitivity due to the 
center perturbation, in the final step. We discuss the modification to which this leads in the $n$-loop (section \ref{sec:nloop}). Then, we discuss the modification 
in the calculation of the sensitivity in the next subsection. In the $n$-loop, 
in addition to step 2, we must also subtract from $Q_n$ 
its projection along $\mathcal{F}$, which is approximately the 
center direction, as:
$Q_n \longleftarrow Q_n -\mathcal{F}_n \mathcal{F}_n^T Q_n/\norm{\mathcal{F}_n}^2,$ 
where $\mathcal{F}_n := \mathcal{F}(u_n,\mathcal{S}).$
Similarly, for $v_n$, after step 3,
$v_n \longleftarrow v_n - \mathcal{F}_n \mathcal{F}_n^T v_n/\norm{\mathcal{F}_n}^2.$

Next, we discuss modifications to adjoint shadowing due to the 
center perturbation. As we noted earlier, the tangent and adjoint subspaces 
corresponding to two different Lyapunov exponents are 
perpendicular to each other (see e.g. Appendix B in \cite{angxiu-siam} for a proof; we also 
show numerical results verifying this fact for the time-delayed 
system considered, in Figure \ref{fig:tan_adj}). This orthogonality gives rise 
to the constraint (derived in section 5.4 of \cite{angxiu-siam}): 
$$(1/N)\sum_{n=1}^{N} (v^{\rm sh}_n)^T \mathcal{F}_n = 0.$$
By definition, since $v^{\rm sh}_n = v_n + Q_n a_n,$
this leads to,
\begin{align}
		\sum_{n=1}^N \Big( v_n^T \mathcal{F}_n + a_n^T Q_n^T \mathcal{F}_n\Big)  = 0. 
\end{align}
This condition leads to one more equation (adding one more row to $G$) while 
solving the NILSS problem (section \ref{sec:lss}). The $n$-loop in adjoint shadowing 
need not be modified. 

\subsection{Computation of the sensitivities}
Having obtained the sequences $\left\{a_n\right\}$, $\left\{v_n\right\}$ and 
$\left\{Q_n\right\}$, the 
shadowing perturbation is determined, for $1\leq n\leq N$, as 
\begin{align}
    v^{\rm sh}_n = v_n + Q_n a_n.
\end{align}
With the shadowing perturbation, the sensitivities can 
be computed as though the system were not chaotic (Eq. \ref{eqn:tangentSensitivity} and 
Eq. \ref{eqn:adjointSensitivity}). That is, in adjoint shadowing,
\begin{align}
		d_s\langle J\rangle = \dfrac{1}{N} \sum_{n=1}^N v^{\rm sh}_n \cdot x_{n'+1}.
\end{align}
In the case of tangent shadowing, we add to the sensitivity in Eq. \ref{eqn:tangentSensitivity}, 
the contribution from the center direction, if treating the center direction separately 
as described in section \ref{sec:center},
\begin{align}
    \label{eqn:tangentNILSSSensitivityFinal}
		d_s\langle J\rangle = \dfrac{1}{N} \Big( \sum_{n=1}^N x^*_n\cdot v^{\rm sh}_n + \frac{v^T_n \mathcal{F}_n}{\norm{\mathcal{F}_n}^2} (J_n - \langle J\rangle_N)\Big),
\end{align}
where $J_n$ in the above equation is the objective function at $u_n$. 
In both tangent and adjoint shadowing,
the projections onto the center direction, $v_n^T \mathcal{F}_n$ 
and $Q_n^T \mathcal{F}_n$, are stored 
during the $n$-loop. However, the $n$-loop needs to be modified to account 
for the center direction only while performing tangent shadowing.
\subsection{Automatic differentiation}
As discussed in section \ref{sec:adintro}, we can replace tangent/adjoint solvers with forward/reverse-mode automatic differentiation,  respectively. 
In the $n$-loop (section \ref{sec:nloop}), we can introduce AD to advance $Q_n, v_n$. Hence, AD-shadowing only requires the primal solver to be supplied by the user, 
as opposed to primal, tangent and adjoint solvers. The AD-version of the $n$-loop is shown in \ref{sec:appxAD}.

We remark that for the AD version of the shadowing algorithm, an exploration of various techniques for memory and time-efficiency of AD \cite{boyana, huckelheim}, such as combining primal-tangent/adjoint solver, is needed. These approaches may lead to taking longer timesteps without compromising on accuracy by utilizing the fact that AD is an exact method, which does not increase the numerical error in the perturbations. This more involved approach to AD shadowing, must be numerically investigated for a given problem to determine whether (or not) it leads to a realizable computational advantage (due to AD overheads) in practice.  

\section{The time delayed model for thermoacoustics}
\label{sec:rijke}

Chaotic thermoacoustic oscillations originate from two main physical nonlinearities, which are deterministic. 
First, the heat released by the flame is a nonlinear function of the acoustic perturbations at the flame's base, i.e. the flame saturates nonlinearly~\cite{Dowling1997,Dowling1999}.  
Both experimental investigations \cite{Gotoda2011,Kabiraj2011, Gotoda2012, Kabiraj2012} and numerical studies \cite{Waugh2014,Kashinath2013,Orchini2015a} showed that the nonlinear flame saturation may cause a periodic acoustic oscillation to become chaotic, by either period doubling, or Ruelle-Takens-Newhouse, or intermittency scenarios~\cite{Nair2014,Nair2015}, which are common in fluid dynamic systems~\cite{Eckmann1981,Miles1984,Eckmann1985}. The numerical studies of~\cite{Waugh2014,Kashinath2013,Orchini2015a} showed that the nonlinear flame saturation may generate chaotic acoustic oscillations even in laminar flame models, where the turbulent hydrodynamics is not modelled. We introduce 
a nonlinear time-delayed model of chaotic thermoacoustic instabilities. 
We demonstrate that shadowing obtains useful sensitivities
of this model, in the chaotic regime. We begin by describing the flame duct model of 
combustion in a horizontal Rijke tube \cite{sujith}, 
open to the atmosphere on both ends. The inviscid 
momentum and energy equations are linearized about the mean flow to yield,
\begin{align}
    \label{eqn:rijke-mom}
    \dfrac{\partial u}{\partial t} &+ \dfrac{\partial p}{\partial x} = 0 \\
    \dfrac{\partial p}{\partial t} &+ \dfrac{\partial u}{\partial x} + \zeta p - \dot{q}\: \delta(x - x_f) = 0,
    \label{eqn:rijke-ene}
\end{align}
where, $u(x,t)$ and $p(x,t)$ are the acoustic velocity and pressure at the 
one-dimensional spatial location 
$x$ and at time $t.$ The pointwise heat-release source is $\dot{q} \: \delta(x - x_f),$ where $\delta(x - x_f)$ is the Dirac delta centered at $x_f.$  The constant $\zeta$ is a parameter that models acoustic damping \cite{Landau1944, Landau1987}. We consider a Galerkin modal 
decomposition in a Fourier basis \cite{sujith}
of the acoustic velocity and pressure fields, which transforms Eqs. \ref{eqn:rijke-mom}-\ref{eqn:rijke-ene} into a set of time-delayed coupled oscillators
\begin{align}
\label{eqn:rijke}
& \frac{d\eta_j}{dt} - j\:\pi\:{\theta_j} = 0, \\
\notag
& \frac{d\theta_j}{dt} + j\:\pi\:\eta_j + \zeta_j\:\theta_j + 2 \beta \:\dot{q}(u_f(t-\tau))\: \sin(j\pi x_f)=0, 
\end{align}
where $j = 1,\cdots,d_g$, indicates the index of 
the Galerkin modes, with $d_g$ being their total number,  
\begin{align}
& u_f(t) = \sum_{k=1}^{N_g}\eta_k(t)\:\cos(k\pi x_f), \;{\rm and}\;\\
& \zeta_j = c_1 j^2 + c_2 j^{1/2}. 
\end{align}
In the above system of equations, $\eta_j$ indicate the velocity modes and $\theta_j$
indicate the pressure modes. Modal damping is represented by $\zeta_j$ and $\dot{q}$ is 
the rate of heat release at the flame location $x_f$. The function $\dot{q}(u) = \sqrt{|1.0 + u|} - 1$ is a modified King's law \cite{sujith}, which is non-differentiable at $u = -1.$ 
Along trajectories of this system, we may encounter states corresponding to 
$u_f = -1$, where the Jacobian does 
not exist. We need to avoid this 
non-differentiability in order to compute the tangent/adjoint/AD 
solutions needed for the shadowing algorithms. Thus, we follow 
the approach taken in \cite{Huhn2020_jfm}, wherein the function 
$\dot{q}$ is approximated by a polynomial around the non-differential point. In particular, when $-1.01 \leq u \leq -0.99$, we take 
$$\dot{q}(u) = -1 + 1750\; (1 + u)^2 - 7.5\times 10^6 \;(1+u)^4,$$
where the coefficients have been estimated by regression. The flame velocity $u_f$ affects the pressure field modes
$\theta_j$ after a time delay given by a constant, $\tau$. This models the fact that the disturbances in the flame velocity at the flame base require a finite time to traverse the flame and cause a perturbation in the heat released \cite{Lieuwen2013, Magri2013}.
\subsection{Replacing the time delay with an advection equation}
When solving the above system numerically, if the time delay
is modelled by converting the delay parameter $\tau$ into 
an integer number of timesteps -- that is, $\tau$ is converted 
into a discrete parameter from a continuous one -- the state 
cannot be differentiated with respect to $\tau$.
We resolve this problem in order to ensure that AD/tangent/adjoint solvers are  
applicable, by augmenting
the primal system with an auxiliary linear advection model \cite{Huhn2020_jfm}
\begin{align}
\label{eqn:advection}
\tau\frac{\partial v}{\partial t} + 2 \frac{\partial v}{\partial y} &= 0 ,\quad
-1 \le y\le 1 \\
v(y=-1,t) &= u_f(t).
\end{align}
The exact solution of the above advection equation  
at the right boundary is $v(y=1,t) = u_f(t-\tau)$. The advection solution, $v(y=1,t)$ can 
be used in place of the heat release model, which in turn influences the pressure modes as per
Eq. \ref{eqn:rijke}. Thus, we mathematically make the overall primal system memory-less, by including
the advection subsystem in Eq. \ref{eqn:advection} in the primal system (Eq. \ref{eqn:rijke}). 
The above discretized equation retains the chain rule dependence on $\tau$ and, hence,
can be differentiated with respect to $\tau$ through AD. We use
a Chebyshev spectral collocation method 
(\cite{Trefethen2005} Ch. 6) to solve the advection equation.
The additional cost per timestep incurred due to adding this advection subsystem 
(i.e., adding Eq. \ref{eqn:advection} to the primal set of ODEs 
in Eq. \ref{eqn:rijke}),
depends on the spatial scale (in $y$) of the numerical discretization of 
the advection system, i.e., the number of Chebyshev collocation points. 
We seek to minimize the number of Chebyshev points in order keep the overall 
dimension of the system as small as possible. The timestep size can also be 
commensurately increased, ensuring the CFL condition, on decreasing the spatial 
resolution, which is cost-effective for computing long-time averages. 
With these considerations, we choose $d_c = 10$ Chebyshev points in the interval $-1\leq y\leq 1$ and the timestep of the primal system is chosen to be $\tau/(2 N_c).$ Choosing $d_g = 10$, the primal system of dimension $d = 2 d_g + d_c = 30$ is time-evolved by integrating Eq. \ref{eqn:rijke} and Eq. \ref{eqn:advection} using the Tsitouras Runge-Kutta time-integrator (Tsit5()) 
offered by the Julia package OrdinaryDiffEq \cite{ordinaryDiffEq, diffeqflux}. We remark that this auxiliary equation approach can be used in any general system
with a constant time-delay in order to maintain its differentiability through
AD with respect to the delay parameter.

\subsection{Types of solutions over a range of the heat-release parameter}
We fix the damping coefficients at $c_1 = 0.06$, $c_2 = 0.01$, the delay parameter 
at $\tau = 0.2$ \cite{Magri2019_amr}, and numerically solve the primal system. 
In this section, we study the effect 
of the heat-release parameter $\beta$
on the \emph{type} of primal solution observed. Here, ``type'' refers to the three 
different possibilities for the asymptotic behavior of a nonlinear dissipative dynamical system, apart from convergence to a fixed point: convergence to a (a) limit cycle (periodic behavior), 
(b) quasiperiodic attractor and (c) chaotic attractor. These different 
regimes are all observed upon varying 
the parameter $\beta$ from about 2 to 9; 
smaller values of $\beta \lesssim 0.8$ lead to 
a fixed point solution. We show these regimes as a function of 
$\beta$ in the bifurcation diagram of Figure \ref{fig:bif}. 

\subsubsection{Ergodic average of acoustic energy}
\begin{figure}
	\includegraphics[width=\textwidth]{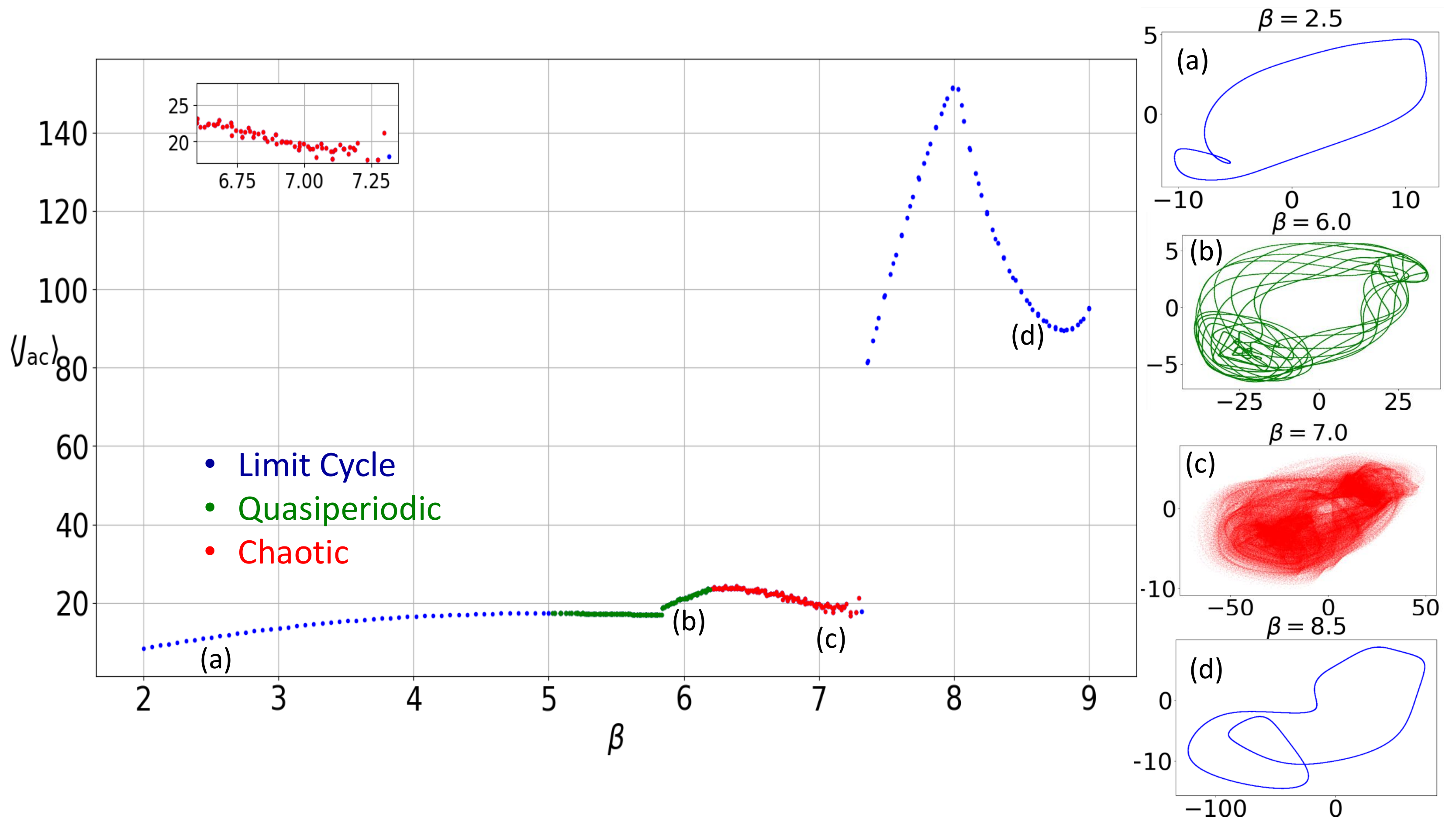}
	\caption{Bifurcation diagram. The values of $\langle J_{\rm ac}\rangle$, the time-averaged acoustic energy, are color-coded according to the type of solution: periodic (blue), quasiperiodic (green) and chaotic (red). 
	The attractors represented on the $u_f$-$\dot{q}$ plane, at four different $\beta$ values -- $\beta = 2.5, 6.0, 7.0, 8.5$ (from top to bottom) -- are shown to the right. The attractors shown on the $u_f$-$\dot{q}$ plane, are also color-coded according to their type. On the top-left inset figure is shown a zoomed-in plot of $\langle J_{\rm ac}\rangle $-vs.-$\beta$ at $\beta = 7.0.$}
	\label{fig:bif}
\end{figure}

The bifurcation diagram
in Figure \ref{fig:bif} shows the acoustic energy upon time-averaging over a long time 
window against the heat release parameter $\beta$. 
This quantity is commonly used as an objective function for optimization problems in thermoacoustics, 
and shall therefore be used to demonstrate the discrete shadowing algorithm.
We denote the acoustic energy by $J_{\rm ac}$, and its ergodic/ensemble average 
by $\langle J_{\rm ac}\rangle$, the latter quantity being computed numerically by time-averaging over a 
\emph{long} trajectory. 
The instantaneous acoustic energy is the sum of the acoustic kinetic and potential energies, i.e., it is the Hamiltonian (constant of motion) of the natural acoustic system. 
Using Parseval's theorem, the acoustic energy is related to the Galerkin modes.
That is, the acoustic energy
is defined as 
\begin{align}
\label{eqn:Jac}
    J_{\rm ac} = \dfrac{1}{2}(p^2 + u^2) = \dfrac{1}{4} \sum_{j=1}^{N_g}(\eta_j^2 + \theta_j^2).
\end{align} The length of the averaging window is chosen to be the time taken for the standard error in the empirical mean values to be within 1\% of the mean, when computed in the chaotic regime. This time, the time-average over which approximates the infinite-time/ensemble average of a given function, was determined to be about 200 time units (or 20000 timesteps, with a fixed timestep size of 0.01).

Another common objective function is the Rayleigh index, whose long-term behavior is also the subject of the sensitivity studies in this paper. The Rayleigh index is defined as \cite{luca-rayleigh, Huhn2020_jfm}
\begin{align}
\label{eqn:Jray}
    J_{\rm ray} := p(x_f, t)\: \dot{q}(t) = \dfrac{1}{2}\sum_{j=1}^{N_g} 
\zeta_j \:\theta_j^2.
\end{align} The physical 
significance of both these objective functions is discussed later in section \ref{sec:obj}.
\subsubsection{Limit cycles and quasiperiodicity}
From Figure 
9, it can be seen that $\langle J_{\rm ac}\rangle$ 
increases continuously over a range 
of $\beta$ values from 2 to about 5. In
this range, the primal 
solution is a limit cycle, which has the effect that the values of $\langle J_{\rm ac}\rangle$ (shown as blue dots) appear to be perfectly observed, without any noise. With $d_g = 10$ Galerkin modes and $d_c = 10$ Chebyshev points, the periodic attractor lives in a 30-dimensional space. For visualization, we show 2D phase diagrams on the $u_f$-$\dot{q}$
plane, for the different solution regimes.
The limit cycle phase diagram is shown in blue in the top-right of the
bifurcation diagram. 

When $\beta \gtrsim 5$, the limit cycle transitions into 
quasiperiodic oscillations, which are aperiodic but appear to be \emph{almost} periodic. For example, an iterative process of rotation on a complex unit circle (or more generally, on the surface of a $d$-dimensional torus) by a constant rational angle is periodic, while a rotation by a constant irrational angle, is quasiperiodic. Mathematically, a quasiperiodic solution is distinguished from a periodic solution by the number of zero Lyapunov exponents:  quasiperiodic solutions have more than one while periodic solutions have exactly one zero Lyapunov exponent. We compute the Lyapunov exponents numerically (using a standard algorithm as explained in section \ref{sec:nloop}) in 
order to classify the different types of solutions \cite{Huhn2020_jfm}. 

The quasiperiodic phase diagram, and the $\langle J_{\rm ac}\rangle$ values in this regime are color-coded green in Figure \ref{fig:bif}. Quasiperiodicity occurs in the transition from periodic behavior to chaotic behavior.
Both the periodic and quasiperiodic case are nonlinearly stable, i.e., the nonzero Lyapunov exponents are negative. This means 
that an applied (infinitesimal) perturbation does not grow exponentially (it may have subexponential growth) in either case. 

\subsubsection{The chaotic regime}
When $6.4\lesssim \beta \lesssim 7.3,$ 
the solutions exhibit at least one positive Lyapunov exponent: this is the chaotic regime. In the bifurcation diagram (Figure \ref{fig:bif}), this regime is shown in red. The phase portrait on the $u_f$-$\dot{q}$ plane shows, as expected, a fractal attractor. The values of $\langle J_{\rm ac}\rangle$ also appear to be erratic, revealing the presence of statistical noise due to a finite time-averaging window. In Figure \ref{fig:LEs}, we show the first 20 Lyapunov exponents at $\beta=7.$ The value of the first exponent is about 0.2. As mentioned in section \ref{sec:decomposition}, the time-derivative of the ODE
is a center perturbation with a zero Lyapunov exponent. The second Lyapunov exponent is  
is about 0.05 corresponding to this center 
direction (it would converge to zero as the averaging time approaches infinity).

\subsection{Acoustic energy and Rayleigh criterion as objective functions}
\label{sec:obj}
We analyze the chaotic thermoacoustic oscillation of the primal 
system by studying the sensitivities of the long-time averages $\langle J_{\rm ac}\rangle$ 
and $\langle J_{\rm ray}\rangle$. Before we compute the sensitivities, we motivate our
particular choice of objective functions, among many available candidates for norms 
~\cite{Chu1965,George2012}, semi-norms~\cite{MagriPhD,Blumenthal2016}, and physical 
measurements in this multi-physical system.
\begin{figure}
    \centering
	\includegraphics[width=0.8\textwidth]{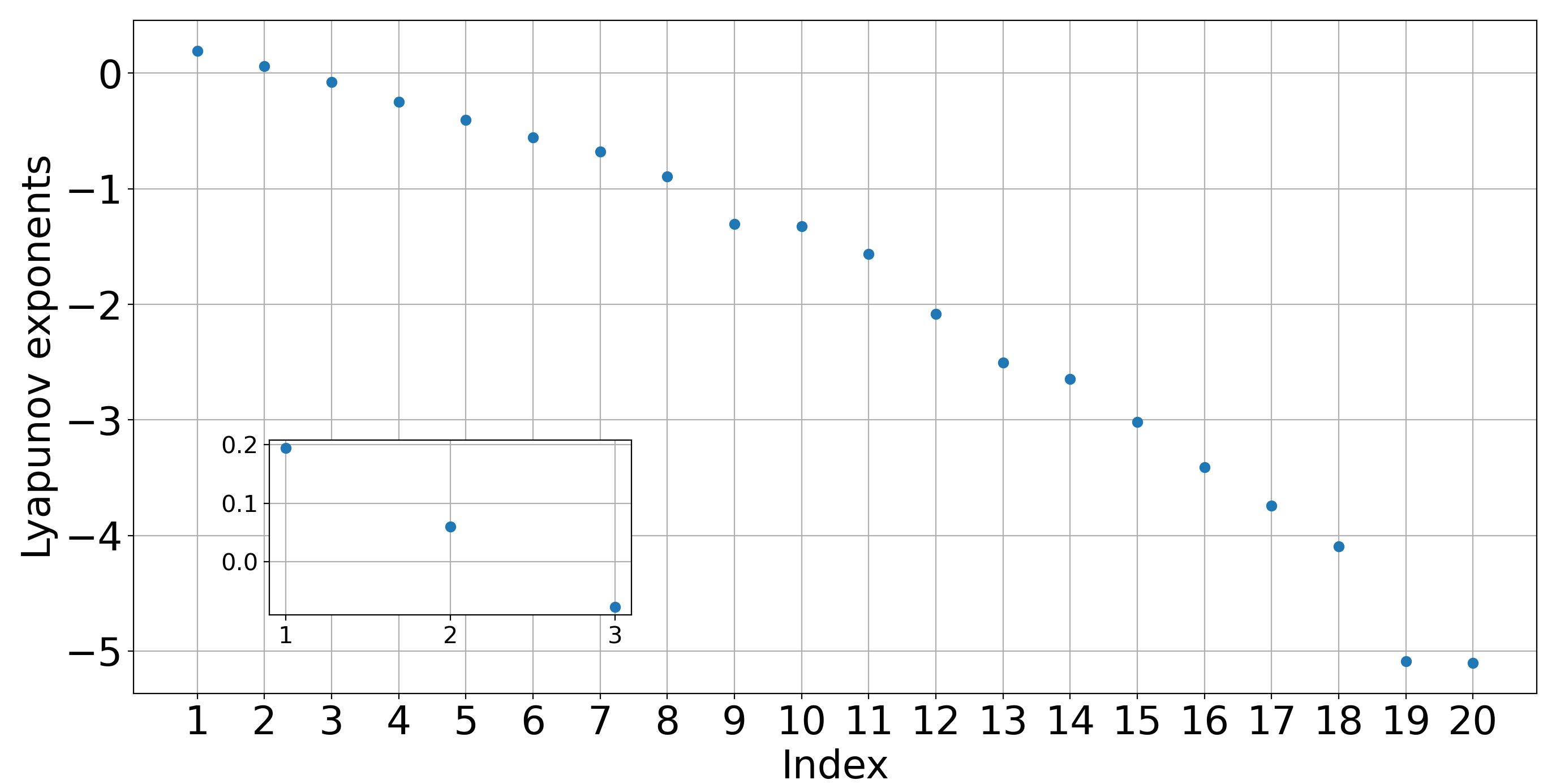}
	\caption{The first 20 LEs at $\beta = 7.0$ and 
		$\tau = 0.2$. The QR factorization is performed 
		every timestep, ie, at a segment length of 
		0.01. The total integration time considered 
		is 200 time units.
	Inset: the first 3 LEs at $\beta = 7.0$ and 
		$\tau = 0.2$. The values obtained are 
		$\lambda_1 \approx 0.19$, $\lambda_2 \approx 0.05$, 
		and $\lambda_3 \approx -0.07$.}
	\label{fig:LEs}
\end{figure}
For thermoacoustic systems with negligible mean flow, which cannot advect flow inhomogeneities like entropy spots, the acoustic energy and Rayleigh criterion are two suitable quantities of interest \cite{luca-rayleigh}. 

%
Since $J_{\rm ac}$ is (half) the Euclidean semi-norm of the thermoacoustic system, we are interested in calculating 
the sensitivity of its time average, $\langle J_{\rm ac} \rangle$, in the interests of reducing the amplitude of chaotic oscillations.
The Rayleigh index can be derived by (i) multiplying the acoustic momentum equation \eqref{eqn:rijke-mom} by $u$; (ii) multiplying the acoustic energy equation \eqref{eqn:rijke-ene} by $p$; (iii) adding them up; and (iv) integrating in the space domain. This procedure yields an equation for the evolution of the acoustic energy  
\begin{equation}\label{eqn:deacdt}
    \frac{dJ_{\rm ac}}{dt} = - \int_0^1 \zeta p^2 \, dx + p_f \dot q,
\end{equation}
where $p_f(t) := p(x_f,t)$ 
is the pressure at the heat source.
Defining the Rayleigh index as 
\begin{align}
J_{\rm ray} := p_f \dot q(t), 
\end{align} 
upon numerical discretization, we obtain    
\begin{align}
J_{\rm ray} = - \dot q(v(y=1,t)) \sum_{j=1}^{N_g} \theta_j(t) \sin(j \pi x_f).  
\end{align}
The Rayleigh index is an important cost functional that determines the stability of acoustic oscillations fed by a heat source. Physically, Eq. \ref{eqn:deacdt} states that the acoustic energy grows in time when the pressure at the heat source is sufficiently in phase with the heat release rate to exceed damping mechanisms. The acoustic energy grows up to nonlinear saturation, after which the self-sustained acoustic oscillation persists. 
This mechanism is commonly studied through the Rayleigh criterion~\cite{Rayleigh1878} for the production of thermoacoustic oscillations. 
In chaotic oscillations, we are interested in calculating the sensitivity of the time-averaged Rayleigh index, $\langle J_{\rm ray} \rangle$.
Applying the infinite time average to Eq. \ref{eqn:deacdt} \cite{Huhn2020_jfm}, 
\begin{align}
   0&= \left\langle \frac{dJ_{\rm ac}}{dt} \right\rangle + \left\langle \int_0^1 \zeta p^2 \, dx + p_f \dot q \right\rangle \nonumber \\
    &= \lim_{T \rightarrow \infty} \frac{1}{T} \int_0^T \frac{dJ_{\rm ac}}{dt} dt + \left\langle  \int_0^1 \zeta p^2 \, dx \right\rangle - \left\langle p_f \dot q \right\rangle \nonumber \\
    \label{eqn:ergAvgOfEac}
    &= \lim_{T \rightarrow \infty} \frac{J_{\rm ac}(T) - J_{\rm ac}(0)}{T} +\left\langle \int_0^1 \zeta p^2 \, dx \right\rangle - \left\langle p_f \dot q \right\rangle. 
\end{align}
Considering that the acoustic energy is a bounded quantity on a strange attractor, the first term of the above equation is 0. 
Hence, Eq. \ref{eqn:ergAvgOfEac} physically means that the damping mechanism exactly balances the acoustic source at regime, i.e., 
\begin{equation}
    \label{eq:rayleigh}
   \langle J_{\rm ray} \rangle = \langle p_f \dot q \rangle = \left\langle {\int_0^1 \zeta p^2 \, dx }  \right\rangle.
\end{equation}
\begin{figure}
    \centering
    	\includegraphics[width=0.8\textwidth]{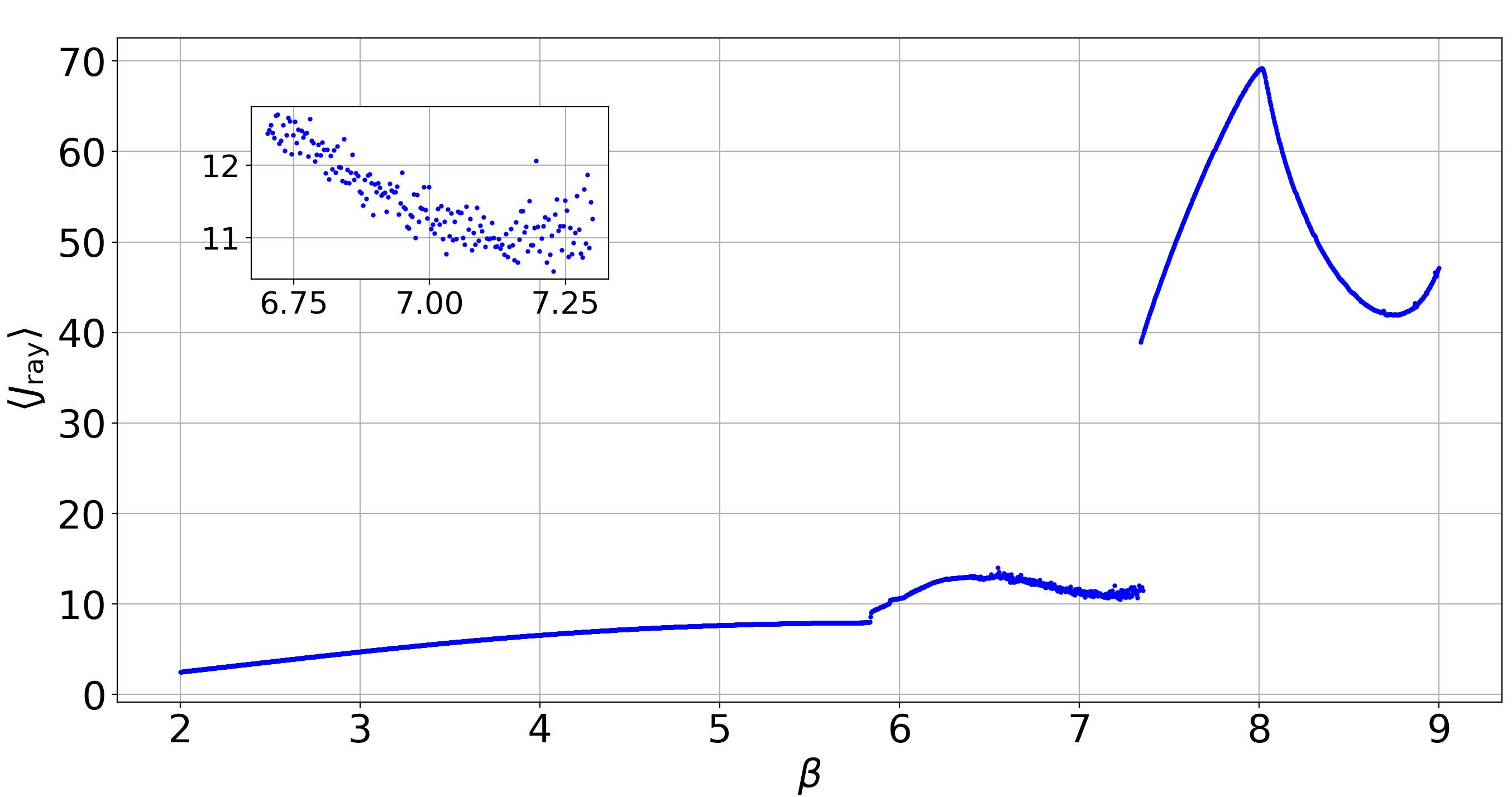}
    \caption{Ergodic average of the Rayleigh index, $\langle J\rangle_{\rm ray}.$ Inset: zoomed-in plot of $\langle J\rangle_{\rm ray}$-vs.-$\beta$ at $\beta = 7.0$ in the chaotic regime. }
    \label{fig:Jray}
\end{figure}
Thus, the time-averaged Rayleigh index can be expressed either from the heat-source contribution or the dissipation term. 
From a computational point of view, the calculation of the sensitivity of $\langle p_f \dot q\rangle$ is difficult because the chaotic modulation, which is imposed exactly at $x=x_f$, makes $\langle p_f \dot q\rangle$ erratic.
To overcome this computational problem, we recommend using $\left\langle {\int_0^1 \zeta p^2 \, dx }  \right\rangle$ (bearing in mind the equality Eq. \ref{eq:rayleigh}), which numerically behaves regularly because it is an integral quantity \cite{Huhn2020_jfm}. Using the Galerkin modal decomposition of the pressure, this integral becomes Eq. \ref{eqn:Jray}, which was used earlier to define the Rayleigh index. In Figure \ref{fig:LEs}, on the right hand side, we plot the ergodic average of the Rayleigh index, $\langle J_{\rm ray}\rangle$ as a function of heat release $\beta.$ The behavior of $\langle J_{\rm ray}\rangle$ is consistent with that of $\langle J_{\rm ac}\rangle,$ with irregular values in the chaotic regime ($6.4 \leq \beta \leq 7.3$), and a sharp increase in the chaotic-to-periodic transition.  
Note that the cost functional $\langle J_{\rm ray} \rangle$ is not directly proportional to the norm of the state, unlike the acoustic energy. The sensitivity and optimization framework we propose can tackle general cost functionals.

\subsection{Lyapunov vectors in the chaotic regime}
\begin{figure}
\centering
	\includegraphics[width=0.4\textwidth]{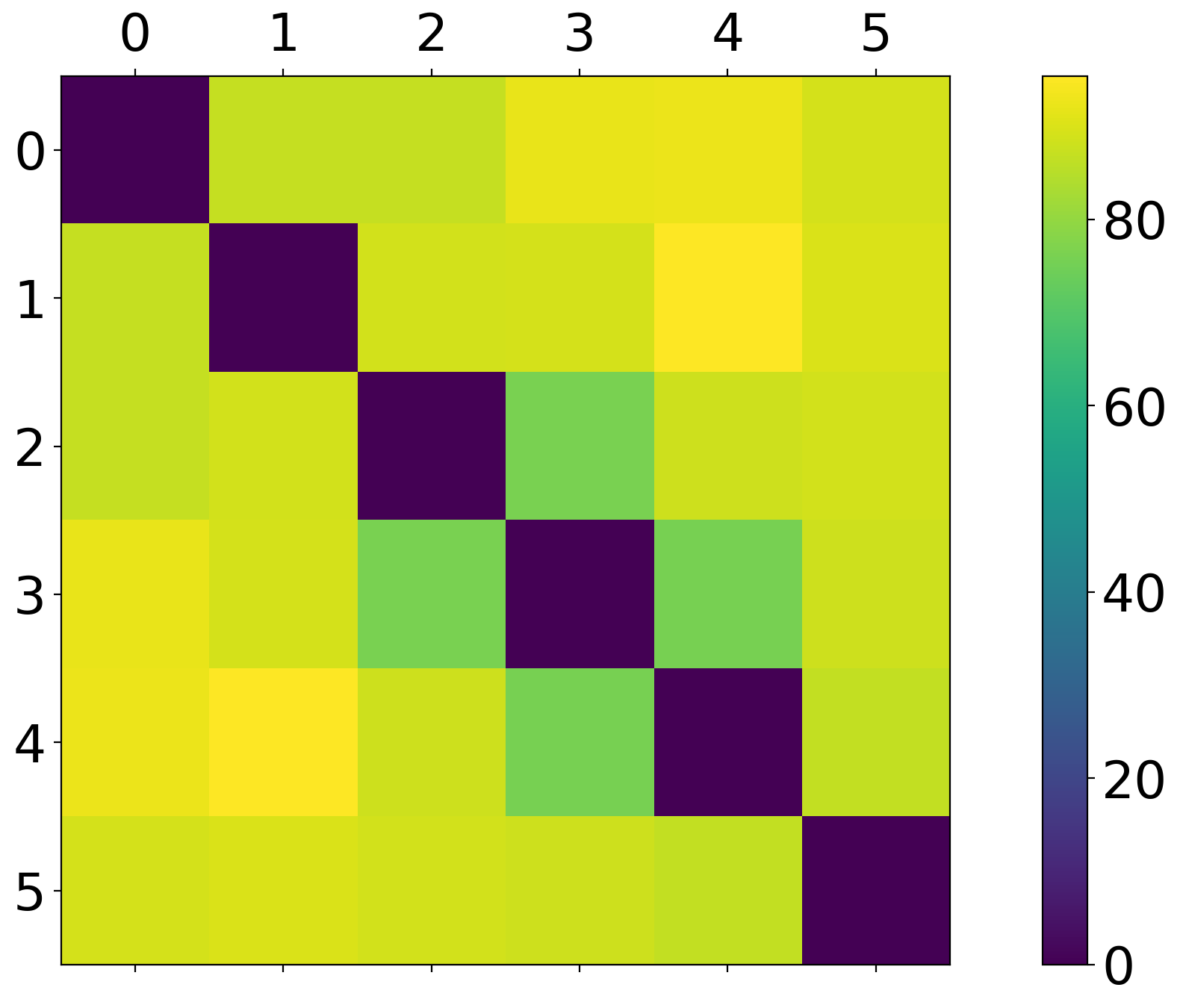}
	\includegraphics[width=0.4\textwidth]{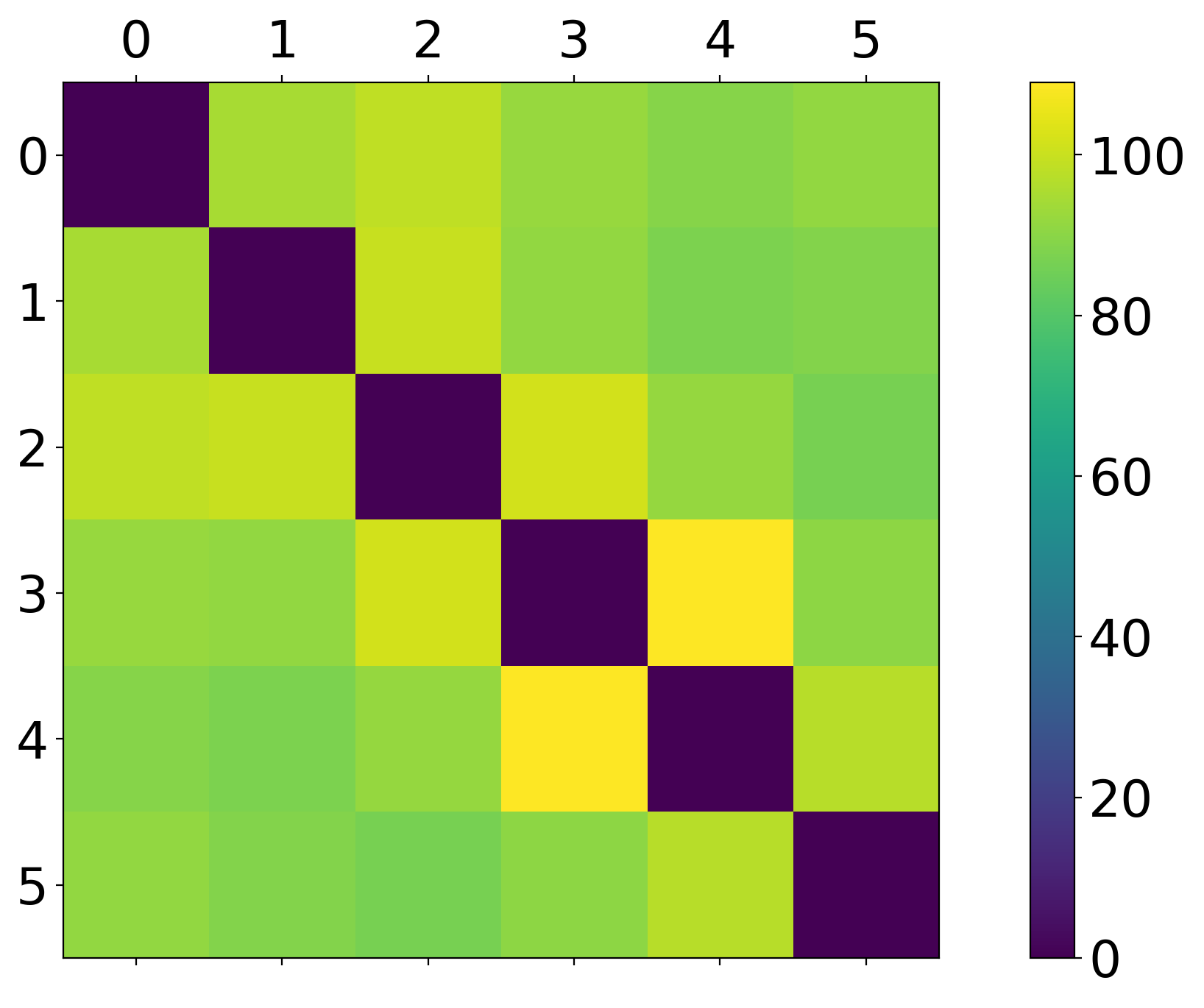}
	\caption{Ergodic average of the
	angles between the first 6 different adjoint CLVs (right) and tangent CLVs (left). The averaging window was set at 250 time units. The colorbar shows the angle in degrees.}
	\label{fig:mean_angles}
\end{figure}

\begin{figure}
\centering
	\includegraphics[width=0.5\textwidth]{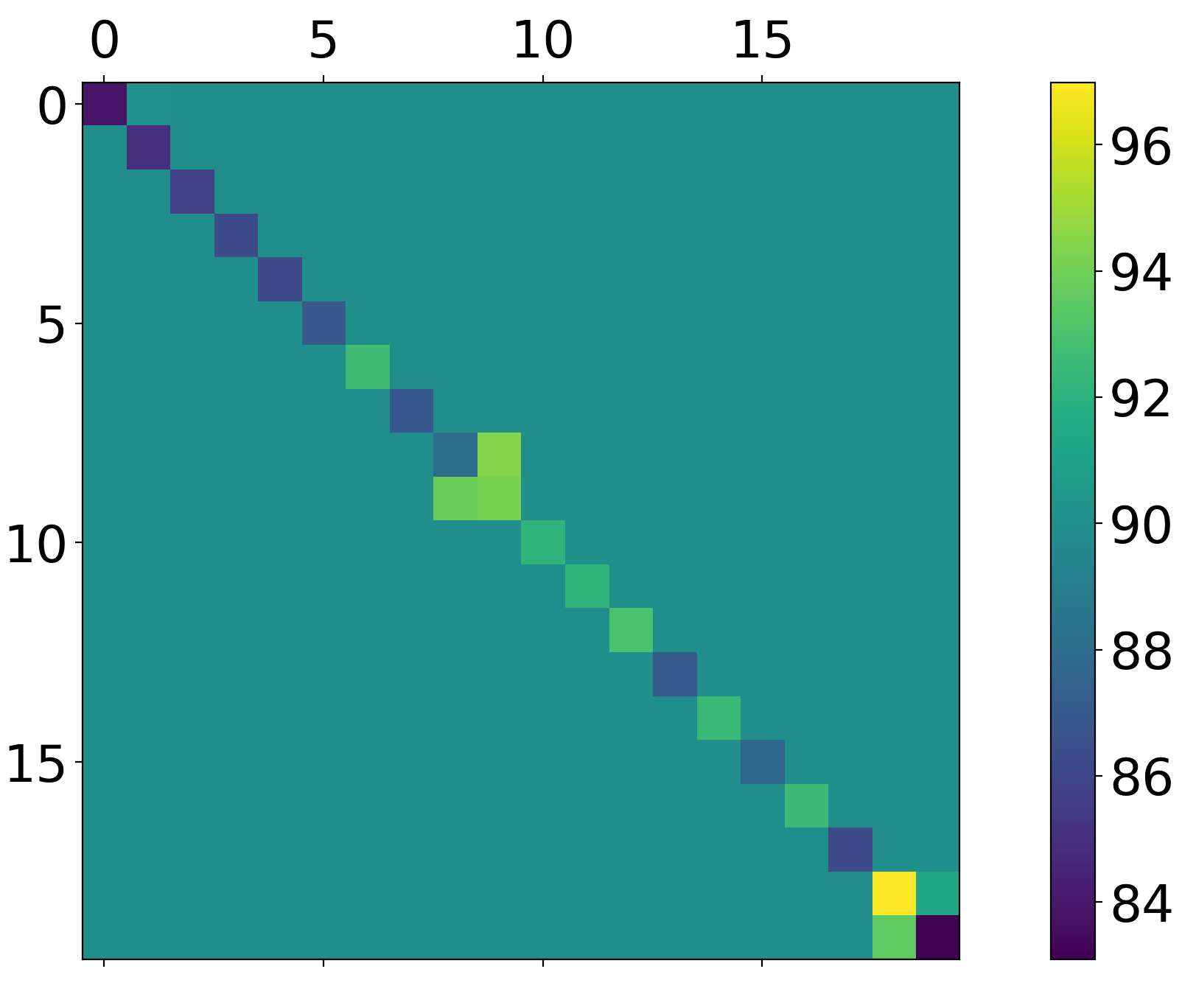}
	\caption{Angles between pairs of adjoint and tangent CLVs when averaged over 
	250 time units. As expected, each tangent CLV is orthogonal to every adjoint CLV except those with the same LE.}
	\label{fig:tan_adj}
\end{figure}

We treat the numerical solution of the system of ODEs as the map $f$ between consecutive timesteps. The parameters $\mathcal{S}$, as per our notation in section \ref{sec:nilss}, is set to 
$[\beta, \tau]^T.$
Corresponding to the LEs (shown in Figure \ref{fig:LEs}), we also compute the  tangent Covariant Lyapunov Vectors (CLVs) using Ginelli {\emph et al.}'s algorithm \cite{Ginelli2007}. The CLVs are tangent/adjoint vectors whose asymptotic exponential growth or decay rates are exactly equal to the LEs; further, they are covariant in the sense that a homogeneous tangent/adjoint solution (Eq. \ref{eqn:homotangent} and Eq. \ref{eqn:homoAdjoint}) starting from a CLV always lies in the span of the same CLV field. The reader is referred to \cite{Kuptsov2012} for the properties of CLVs. In this paper, the CLVs are normalized. The unstable (stable) tangent CLVs form a basis of (non-orthonormal) unit vectors for $E^u$ ($E^s$); the span of the unstable (stable) adjoint CLVs is $(E^s)^\perp$ ($(E^u)^\perp).$ Since the center subspace is one-dimensional, the normalized time-derivative $\mathcal{F}/\norm{\mathcal{F}}$ is the center tangent CLV field. The Ginelli algorithm can also compute the adjoint CLVs. For this, we use the Jacobian transpose trajectory, in place of the Jacobian trajectory used to compute the tangent CLVs, and we also reverse time (i.e, the forward/backward phase of Ginelli's algorithm is carried out backward/forward). We compute the tangent and adjoint CLVs at $\beta = 7$. In Figure \ref{fig:mean_angles}, we show the 
angles between each pair among the first 6 tangent CLVs, and each pair of adjoint CLVs, on the right. The angles are averaged over 250 time units. In a hyperbolic system, the angles between every pair of CLVs (corresponding to different LEs) are uniformly (in phase space) bounded away from 0. Although not a rigorous test for hyperbolicity, the results of Figure \ref{fig:mean_angles} indicate that, at least on average, the first 6 CLVs, both tangent and adjoint, do not show tangencies. Over the time window of calculation, the minimum angle observed between any dissimilar pair was about 4 degrees. This indicates that the system is likely uniformly hyperbolic.  Biorthogonality is numerically verified in Figure \ref{fig:tan_adj}. Except along the diagonals, which contain the mean angles between a tangent and adjoint CLV corresponding to the same LE, the mean angles are all about 90 degrees, as expected.

\section{Suppression of a nonlinear oscillation by gradient-based optimization}
\label{sec:optim}
The thermoacoustic model under investigation displays chaotic behavior in the region $6.4\leq \beta \leq 7.3.$ As illustrated in Figure \ref{fig:perts}, conventional methods to compute the sensitivities of the long-time behavior of this model, in this chaotic regime, fail to produce meaningful sensitivities. In this section, we use the shadowing algorithm presented in section \ref{sec:nilss} to enable the computation of these sensitivities. Our goal is to illustrate the potential of the algorithm for practical sensitivity-based optimization and parameter estimation in the regime of chaotic acoustics. 

Given the small number of parameters and objective functions, it is possible to compare the sensitivities computed through the algorithm with the slopes obtained from the bifurcation diagrams in Figures \ref{fig:bif} and \ref{fig:LEs} (right). These comparisons validate the results of our algorithm, which is one of the goals in this section. We demonstrate the usefulness of the computed sensitivities by using them in a gradient descent algorithm to minimize the ergodic average of the acoustic energy. This simple optimization procedure can be used for heat release parameter selection. The NILSS algorithm and its discrete AD variant presented here thus introduce sensitivity-based optimization and parameter estimation to the chaotic regime, more generally in hyperbolic systems with constant time-delays, extending the work of Huhn and Magri \cite{Huhn2020_jfm}.
\subsection{Shadowing sensitivities of the acoustic energy and Rayleigh criterion}
\label{sec:sens}
We apply the shadowing algorithm from section \ref{sec:nilss} to compute both tangent and adjoint shadowing sensitivities. We compute the tangent shadowing direction once to estimate $(d\langle J_{\rm ac}\rangle/d\beta)$ and 
$(d\langle J_{\rm ray}\rangle/d\beta).$ Similarly, we compute an adjoint shadowing direction once, to calculate both the sensitivities $d\langle J_{\rm ac}\rangle/d\beta$ and $d\langle J_{\rm ac}\rangle/d\tau$. Next we define the inputs to the shadowing algorithms. A primal orbit $\left\{ u_n\right\}$ is a sequence of 
30-dimensional solution vectors obtained by time-integrating the primal system 
(Eq. \ref{eqn:rijke} and Eq. \ref{eqn:advection}). The map $f$ is the Tsitouras Runge-Kutta time-integrator that advances a solution state by one timestep. The timestep size is fixed at 0.01. For tangent shadowing, the input $b_n$ is set to $x_n = (\partial f/\partial \beta)(u_{n-1}, \mathcal{S}).$  In the adjoint shadowing algorithm, we arrange the input sequence $\left\{b_n\right\}$ so that $b_{n'+1}$ is set to $x^*_n = DJ_{\rm ac}(u_n).$ (i.e., we pass $x^*_n$ in time-reversed order to the adjoint shadowing algorithm). We use the AD package \verb+Zygote.jl+ \cite{diffeqflux} to compute the sequence $\left\{ b_n \right\}$ for tangent shadowing, through AD. The Jacobian matrix $A_n := Df(u_n)$ and its transpose, needed for the tangent and adjoint algorithms respectively, are computed by using finite difference. Instead, if using AD, each $A_n$ must be computed row-by-row (since \verb+Zygote.jl+ does not support vector-valued outputs), leading to a much larger computation time compared to using finite difference. Moreover, the shadowing algorithms do not need the Jacobian to be computed exactly. In the AD version of tangent and adjoint shadowing, 
we do not need to compute $A_n$ and $b_n$; the tangent/adjoint perturbations needed in the $n$-loop (section \ref{sec:nloop}) are directly computed using AD as shown in \ref{sec:appxAD}. The input to the AD version of both tangent and adjoint shadowing are the functions that perform primal time-integration and compute the objective functions, given the primal state. In order to ensure that $u_0$ is a point on the chaotic attractor, we evolve the system for a time of 10000 time 
units, starting from a random 30 dimensional vector. 

\begin{figure}
	\includegraphics[width=0.5\textwidth]{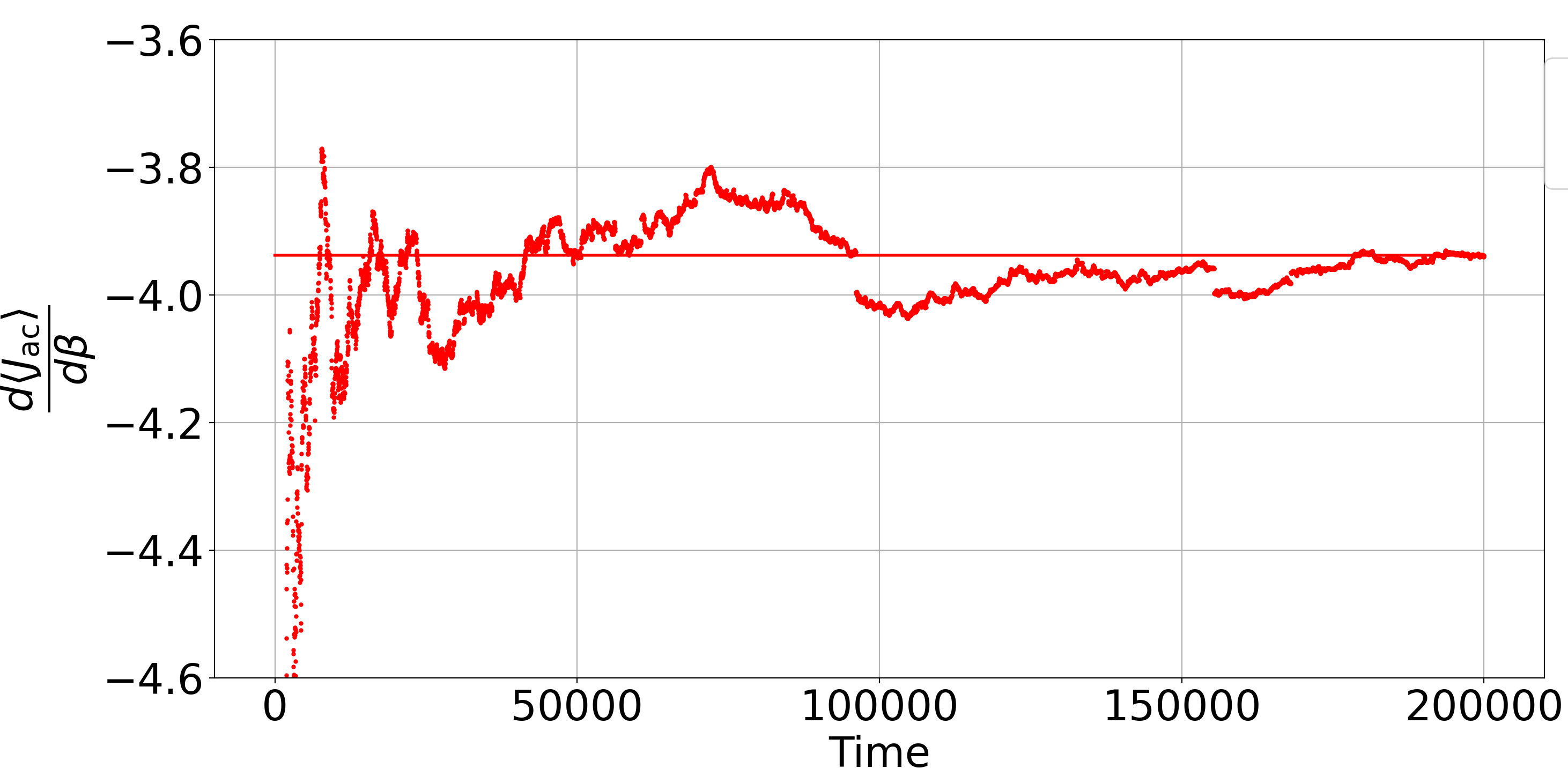}
	\includegraphics[width=0.5\textwidth]{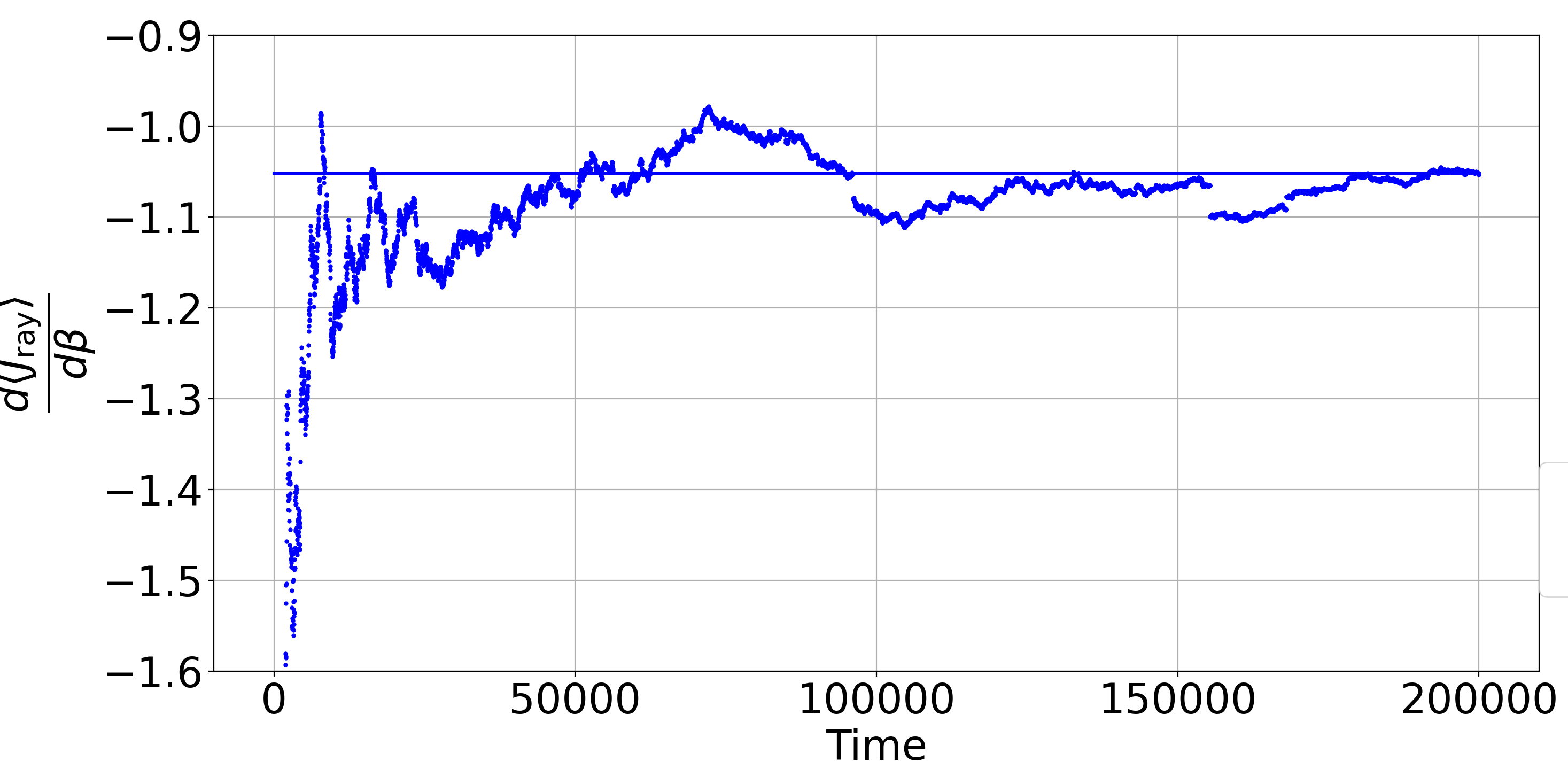}
	\caption{Sensitivities of $\langle J_{\rm ac}\rangle$ (red) and 
	$\langle J_{\rm ray}\rangle$ (blue) computed through tangent shadowing. The sensitivities are obtained by a cumulative average over sensitivities each of which is computed over 20 time units.}
	\label{fig:dJ_dbeta}
\end{figure}

\begin{figure}
	\includegraphics[width=0.5\textwidth]{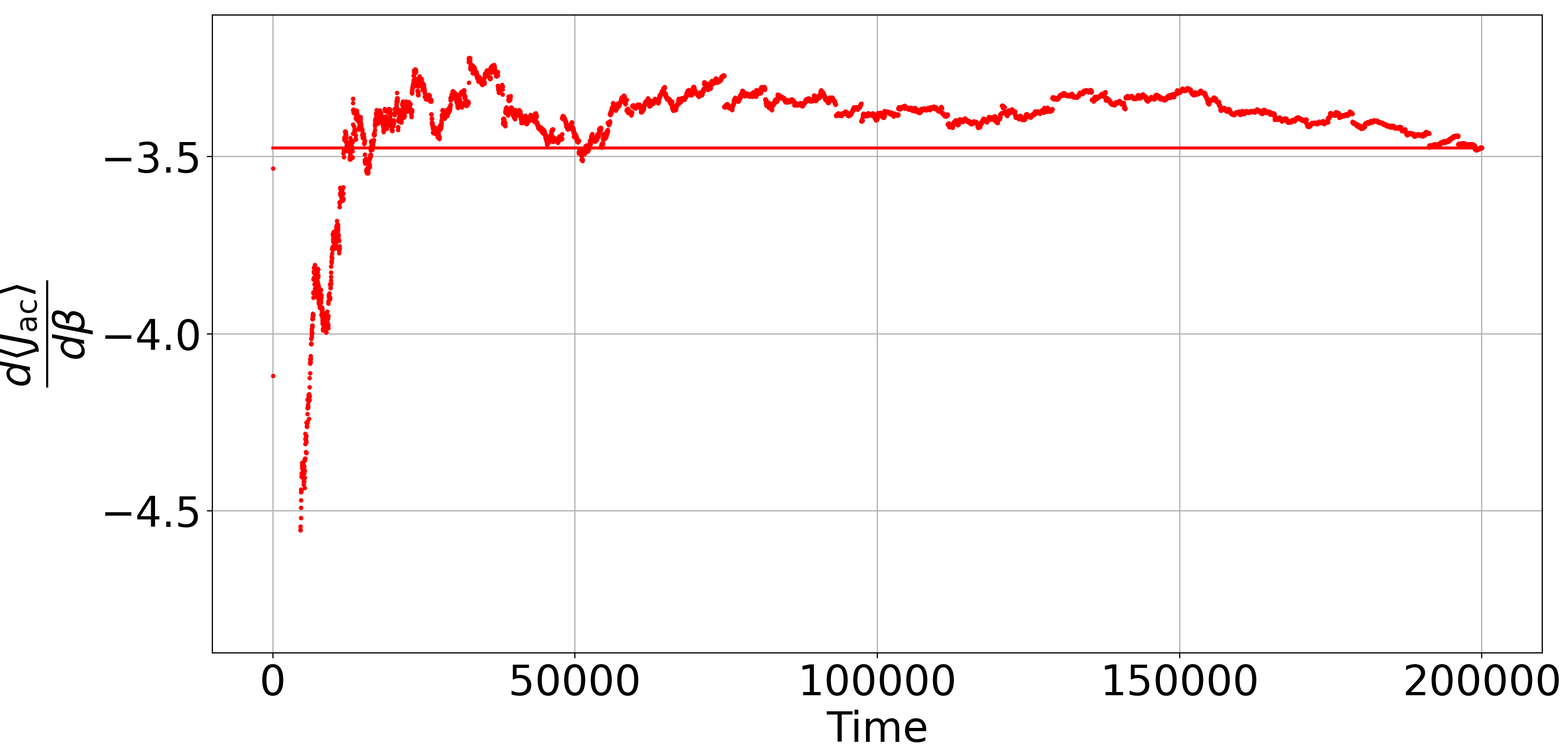}
	\includegraphics[width=0.5\textwidth]{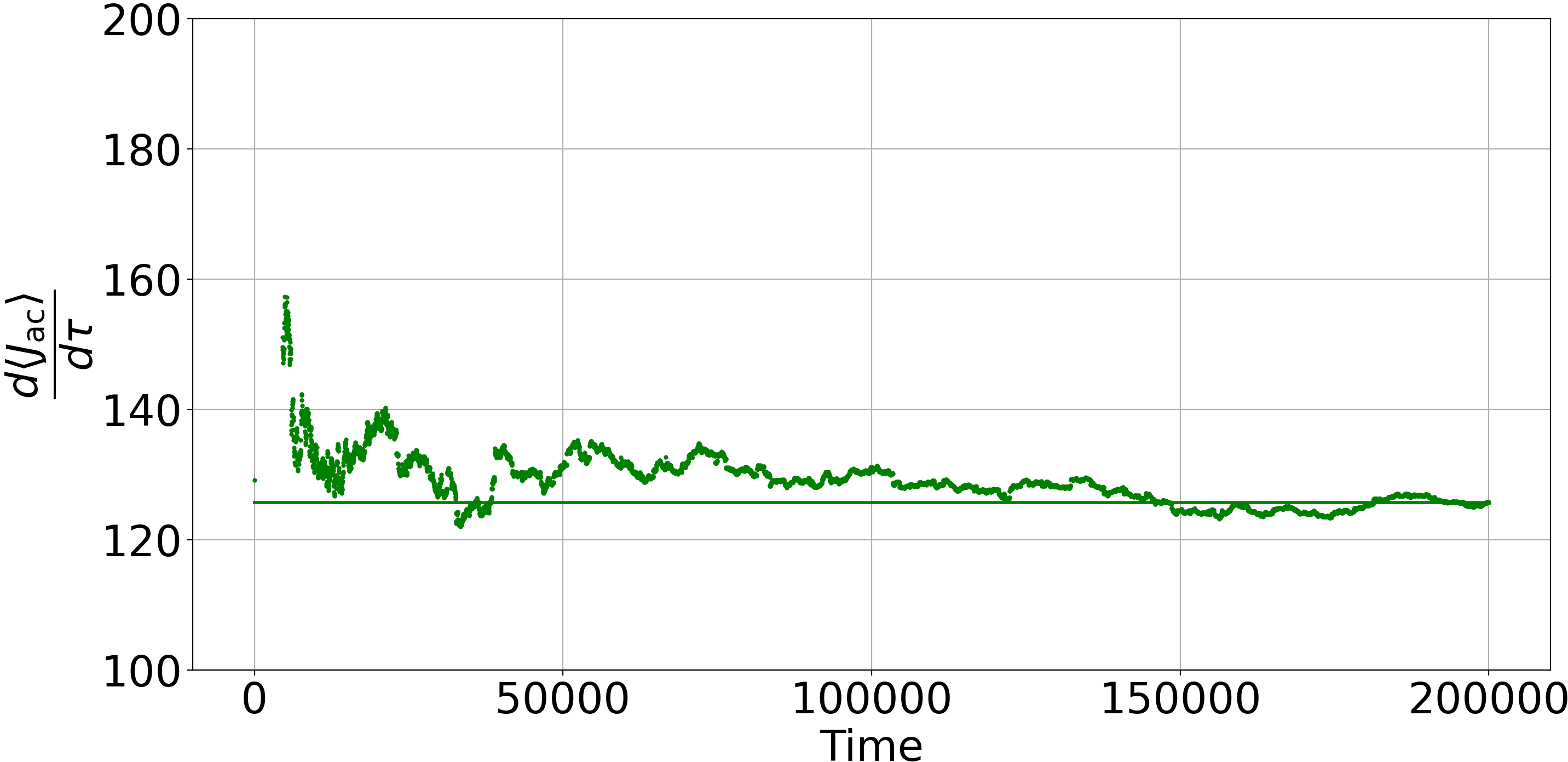}
	\caption{Sensitivities of $\langle J_{\rm ac}\rangle$ with respect to $\tau$ (green, right) and 
	$\langle J_{\rm ac}\rangle$ (red, left), wrt $\beta$, computed through adjoint shadowing at 
	$\beta = 6.9, \tau = 0.2$. The recorded sensitivities 
	were obtained by a cumulative average over sensitivities each computed over a time length of 20 units.}
	\label{fig:dJac_ds}
\end{figure}

We set $d_u = 2$ for both algorithms. Although we could 
theoretically have used $d_u = 1$ in the tangent algorithm, setting $d_u = 2$ leads to 
better approximations of the shadowing direction via the least squares problem (section 
\ref{sec:lss}). Each minimization problem is solved over a time duration of 20 time units, which is about 4 Lyapunov times. To obtain the sensitivity of a 
long-time average, a sample mean of these intermediate-time sensitivities is taken. A 
cumulative mean of the sensitivities converges as the length of time (number of samples) 
increases (\ref{fig:dJ_dbeta} and 
\ref{fig:dJac_ds}). In Figure \ref{fig:dJ_dbeta}, we show the sensitivities of the 
long-time averaged acoustic energy (left) and Rayleigh index (right), with respect to 
$\beta$, computed using sample averages of the tangent shadowing sensitivities.
In Figure \ref{fig:dJac_ds}, we show the sensitivities of the 
long-time averaged acoustic energy with respect to $\beta$ (left)
and $\tau$ (right), computed using sample averages of the adjoint shadowing sensitivities. The mean 
up to a time of 200,000 (i.e., calculated using 10,000 shadowing sensitivities each over a time of 20 units) is shown as a solid line. The mean values shown in both plots in Figure \ref{fig:dJ_dbeta} compare 
well against the corresponding slopes from the inset plots of Figures \ref{fig:bif} and \ref{fig:Jray}. The mean value of the sensitivity with respect to $\beta$, about 
-3.9, also agrees well, with the same sensitivity computed using the adjoint 
algorithm, around -3.5, which is shown in Figure \ref{fig:dJac_ds} (left). Since all 
three quantities, the shadowing sensitivities from the two algorithms as well as a 
direct reading of the slope from the $\langle J_{\rm ac}\rangle$-vs-$\beta$ plots, suffer from statistical noise due to a finite computation window, we do not expect exact agreement. 
Nevertheless, we note that both tangent and adjoint sensitivities are within 12\% of the 
slope estimate of -4 obtained from by approximating $\langle J_{\rm 
ac}\rangle$-vs-$\beta$ as a line in Figure \ref{fig:bif} (inset). To ensure its 
correctness, the program that implements the two shadowing algorithms is validated by computing sensitivities on the classical model of a chaotic ODE, the Lorenz'63 system, to verify the computed sensitivities against the values available for this system in 
the literature. For details on this validation and on the replication of Figures 
\ref{fig:dJac_ds} and \ref{fig:dJ_dbeta}, see \ref{sec:appxLorenzSens}. The 
absolute values of the sensitivities $d\langle J_{\rm ray}\rangle/d\beta$ from Figure 
\ref{fig:dJ_dbeta} are consistent with our expectation from Figure \ref{fig:Jray} that the ergodic average $\langle J_{\rm ray}\rangle$ decreases with $\beta$, around $\beta = 6.9$, but not as rapidly as $\langle J_{\rm ac}\rangle$. 
From Figure \ref{fig:dJac_ds} (right), we see that the ergodic average of the acoustic energy is highly sensitive (when compared to changes in $\beta$) to small perturbations in the time delay parameter $\tau.$ The shadowing algorithm shows convergence, when the same derivative is computationally prohibitive to obtain accurately with ensemble sensitivity calculations \cite{Eyink2004, nisha-ens}.  

\subsection{Minimization of the acoustic energy using shadowing}
\begin{figure}
\centering
	\includegraphics[width=0.9\textwidth]{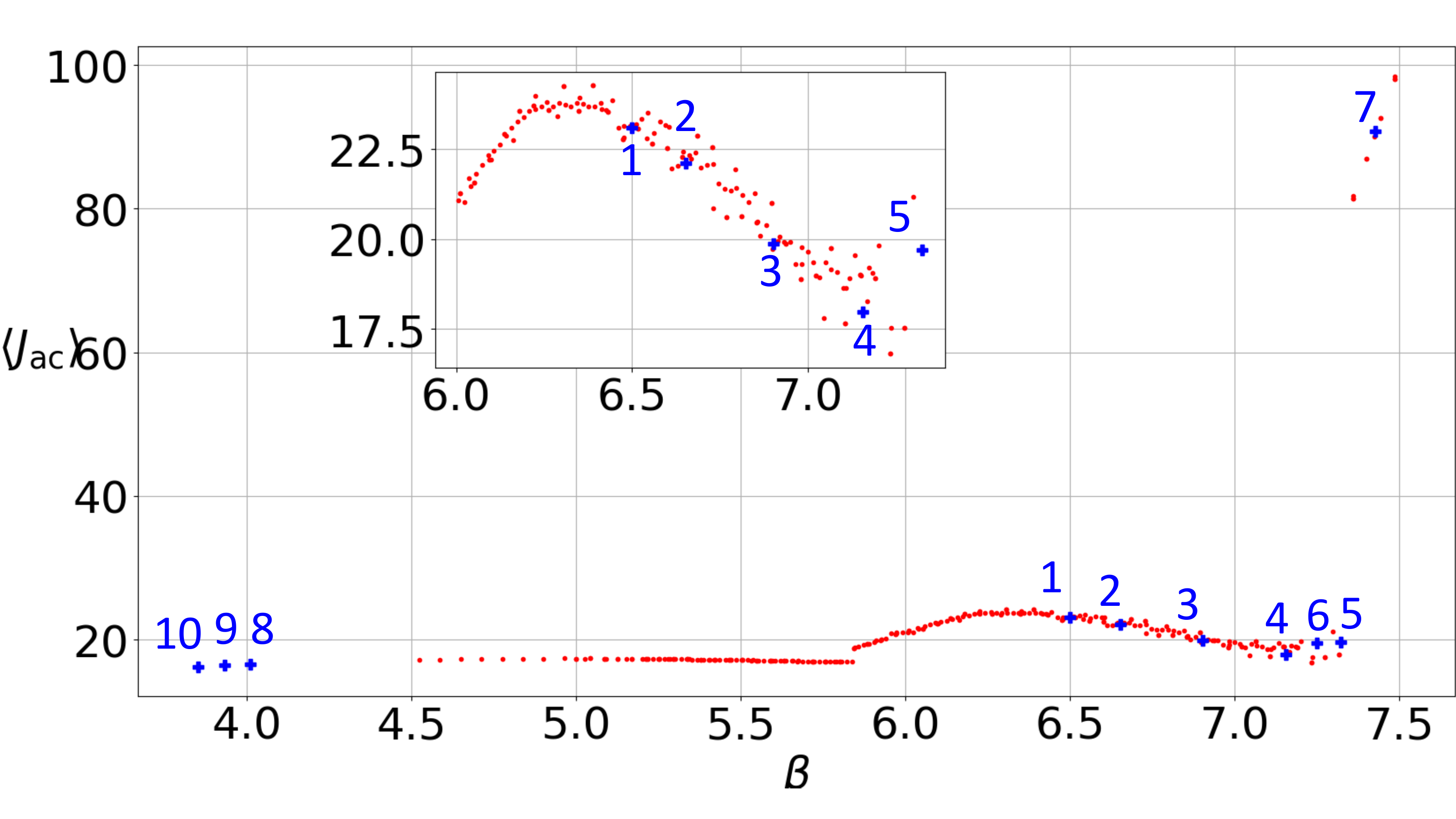}
\caption{The path of optimization starting at $\beta = 6.5$. The points are superimposed on the plot of $\langle J_{\rm ac}\rangle$. In the 
gradient descent algorithm, the step size is taken to be 0.1.} 
	\label{fig:optimization_path}
\end{figure}
We shall demonstrate an application of the sensitivity $d\langle J_{\rm ac}\rangle/d\beta$ that we computed in the previous section, to the problem of optimal design by parameter selection \cite{Huhn2020_jfm}. We use the numerically computed bifurcation diagram (Figure \ref{fig:bif}) as a qualitative check of our parameter selection procedure. We aim to solve the following optimization problem, 
\begin{equation}
    \label{eq:optim_problem}
    \begin{aligned}
        & \underset{\beta}{\text{minimize}} & & \langle J_{\rm ac} \rangle (\beta) \\
        & \text{subject to} & & {\rm Eq.} \ref{eqn:advection},\; {\rm Eq.} \ref{eqn:rijke}, 
    \end{aligned}
\end{equation}
where $\langle J_{\rm ac}\rangle$ is the ergodic average of the acoustic energy, which is approximated using a long-time average over a time duration of $N$ steps.
%
%
%
The parameter can be updated by a steepest descent method  
\begin{equation}
    \beta_{n+1} = \beta_n - \left . \gamma d_\beta \langle J_{\rm ac} \rangle \right |_{\beta =\beta_n},
\end{equation}
A relaxation factor of $\gamma = 0.1$ is used to enable stable and accurate numerical convergence. At each $\beta$, the shadowing sensitivity is computed over 1000 time units, 
which is a sample average of 50 sensitivities collected from runs of the shadowing algorithm each 
over a time of 20 units. The algorithm stops when the condition 
\begin{equation}
    \label{eq:optim_condition}
    \left\langle J_{\rm ac} (\beta_n)\right\rangle < \epsilon \cdot \left\langle J_{\rm ac} (\beta_0)\right\rangle  
\end{equation}
is met, where $\epsilon = 1\%$. This condition physically signifies that the optimization is successful when the system vibrates around the fixed point.
In Figure \ref{fig:optimization_path}, we show the optimization path (blue points) taken by this procedure, starting from $\beta = 6.5$ We show the points $(\beta, \langle J_{\rm ac}\rangle(\beta))$ numbered in the order in which they are encountered in the optimization procedure.
As shown in the figure, the path leads to a reduction in $\langle J_{\rm ac}\rangle$ by 
exiting the chaotic region into the periodic regions for larger of $\beta,$ and eventually 
into the periodic region at smaller values of $\beta$. We remark that, in this case, direct evaluation of the bifurcation diagram is possible, owing to the 
relatively low dimension of the system and the parameter space. This bifurcation diagram itself suggests optimal paths for acoustic energy reduction. However, in a more general setting, the dimension of the system, 
the objective function and parameter spaces may be such that bifurcation diagrams, at the resolution of parameters required to compute accurate gradients, are computationally infeasible. The demonstration in this section 
indicates that the shadowing sensitivities, which can be computed at a smaller cost relative to the bifurcation diagram computation, can be used instead, for optimization problems. For illustration purposes, we have 
chosen a simple gradient-descent with a fixed relaxation factor. The numerical results in Figure \ref{fig:optimization_path}, however, show that the step size plays an important role in the rapid transition from large $\beta$ in the periodic region (around 7.4) to the periodic region for smaller values (around 4.0). We defer to a future work the effect of the step size, which is beyond the scope of this paper.  

\section{Data assimilation with discrete shadowing}
\label{sec:da}
A common problem whenever we have incomplete, and often noisy, observations together with a model of a physical system, is to estimate a model trajectory that reproduces the observations. In this section, we explain that the variational formulation of this problem is an application of the shadowing algorithm discussed in section \ref{sec:nilss}. Then, we apply the shadowing algorithm to the time-delayed model to illustrate its potential for data assimilation in chaotic solutions. 

In a data assimilation problem, we are given external measurements of an observable $g$ at a sequence of times, denoted $g^{\rm obs}_{n_k},$ $1\leq k \leq N$. In the data assimilation method 4DVar \cite{Traverso2018}, an initial state $u_0$ is sought 
so that the model trajectory at the observation times, $f_{n_k}(u_0) = u_{n_k}$  produces a sequence $g_{n_k} := g(u_{n_k})$ that closely 
matches the observations $g^{\rm obs}_{n_k}$. If a reliable guess for the initial state, known as \emph{background} $u^{\rm bg}_0$, is available, we desire our predicted initial state to be close to the background. The optimal initial state, known as \emph{analysis}, is obtained by minimizing the following cost functional
\begin{align}
\label{eqn:costFunctional}
    \langle J\rangle(u_0) = \frac{1}{2}(u_0 - u^{\rm bg}_0)^T B^{-1} (u_0 - u^{\rm bg}_0) + 
    \frac{1}{2}\sum_{k=1}^{N} (g^{\rm obs}_{n_k} - g_{n_k})^T 
    L^{-1}(g^{\rm obs}_{n_k} - g_{n_k}).
\end{align}
The first term in the cost functional in Eq. \ref{eqn:costFunctional} corresponds to the misfit 
between the predicted initial state and the background, 
where $B$ is a $d\times d$ background error covariance matrix.
The second term is the misfit between the external observations and values 
of the observable generated by the model, using the predicted state as the initial condition. Given an observable space of dimension $l$, the $l\times l$ matrix $L$ is the observation error covariance matrix. 
Without loss of generality, take $B$ and $L$
to be identity matrices (in $\mathbb{R}^{d\times d}$ and $\mathbb{R}^{l\times l}$ 
respectively), which 
results in the following cost functional: 
\begin{align}
    \label{eqn:simplerCostFunctional}
    \langle J\rangle(u_0) = \frac{1}{2}\norm{u_0 - u^{\rm bg}_0}^2 + \frac{1}{2 }\sum_{k=1}^{N} 
    (g^{\rm obs}_{n_k} - g_{n_k})^2.
\end{align}
The cost functional is minimized using a standard optimization procedure, and an analysis state $u_0^*$ is obtained. The reader is referred to \cite{Traverso2019} for data assimilation of nonchaotic states with a similar thermoacoustic model 
considered in the present paper. When the model 
is chaotic, the standard optimization procedure fails since the gradient of the functional with respect to $u_0$ grows exponentially with the time duration of available observations (Figure \ref{fig:perts}). Thus, assimilation for a time window longer than the typically short Lyapunov time, cannot be achieved using standard optimization methods. 

However, several effective strategies have been proposed, particularly in the field of numerical weather prediction, wherein chaotic models are widely used, that are successful over long assimilation windows. The most popular of these include 4DVar-AUS, in which the \emph{analysis increment} -- the discrepancy added to the state during optimization -- is restricted, 
to the nonstable ($E^u\oplus E^c$) subspace, at every timestep  \cite{Trevisan2007}. Another method is projected shadowing-based data assimilation in which the cost functional is minimized using Newton's method in each step of which the analysis trajectory, as a whole, is updated. Additionally, the updates to the analysis trajectory during 
each Newton iteration is carried out only on the nonstable subspace, and this leads to the economy of the method. The updates to the analysis trajectory on the stable subspace is treated using a different method, known as synchronization \cite{de_leeuw_projected_2018}. In this work, we present an alternative approach that uses the shadowing algorithm described in section \ref{sec:nilss}. The goal of this approach is to compute a model orbit that shadows the pseudo-orbit pertaining to the observations. Our method thus joins the class of shadowing-based 
data assimilation methods offering an alternative formulation that indirectly computes the 
shadowing orbit through NILSS.

\subsection{Tangent NILSS for state estimation with full-state observations}

We introduced NILSS \cite{angxiu1} as a method to differentiate long-time averages with respect to parameters in a chaotic system. How is the method relevant to the problem of state and parameter estimation? The answer to this question lies in the fact that a shadowing trajectory can be obtained as a byproduct of the NILSS method, and shadowing trajectories can be used for state estimation. First, we recognize that, in NILSS, the derivative of the long-time average is computed along a shadowing trajectory at a 
perturbed parameter. Secondly, we relate 
the shadowing perturbation sequence, $\left\{v^{\rm sh}_n\right\}$, that 
is computed by tangent NILSS, to the solution of the 
state and parameter estimation problem.

\subsubsection{Shadowing-based interpretation of NILSS}
NILSS and its adjoint versions use the shadowing lemma (see e.g. 18.1.2 of \cite{katok}) by considering the reference trajectory $\left\{u_n\right\}$ as a pseudo-orbit of $f(\cdot, s + \epsilon).$ According to the shadowing lemma, there is a unique orbit of $f(\cdot, s + \epsilon)$ called the shadowing orbit, that is close to the given reference orbit of $f(\cdot, s)$, $u_n$. The tangent \emph{shadowing perturbation} is the sequence of tangent vectors, along this shadowing orbit, that expresses the discrepancy between the shadowing orbit and the pseudo-orbit $u_n, n \in \mathbb{Z}^+$, in the limit $\epsilon \to 0.$ Hence it remains bounded for all time. A close approximation of the shadowing perturbation, $v^{\rm sh}$, is obtained by solving the following least squares problem (\cite{shadowing}, Theorem LSS)
\begin{align}
\notag 
    v^{\rm sh} &= {\rm arg min}_{w \in \mathcal{V}} 
    \sum_{n=0}^{N-1} \norm{w_n}^2 \\
    \label{eqn:inhomogeneousTangentConstraint}
    {\rm s.t.} \;\; &w_{n+1} = D_uf(u_n) \; w_n + x_{n+1}.  
\end{align}
In the above optimization problem, the search space $\mathcal{V}$ is 
restricted to a subset of $(\mathbb{R}^d)^N$ in which each $w_n$  
can be expressed as $w_n = v_n + Q_n a_n$, i) for some $a_n \in \mathbb{R}^{d_u}$, ii) where $v_n$ is the conventional tangent solution, i.e., solution of Eq. \ref{eqn:tangent} with zero initial condition. This leads the NILSS algorithm \cite{angxiu1} to being more efficient than if the search space was set to $(\mathbb{R}^d)^{N}$ \cite{angxiu1}. 

\subsection{Converting the state estimation problem to a parameter estimation problem}
The cost functional in state estimation consists of two parts: the background error 
and the observation error. We may assume that the background orbit is a pseudo-orbit 
of $f(\cdot, s + \epsilon)$, for some $\epsilon$ around zero. Then, the 
shadowing problem (Eq. \ref{eqn:inhomogeneousTangentConstraint}) minimizes the first part of the state estimation cost functional, which is the background error. In order to minimize the
second part, the observation error, we use the observation error as the
objective function in the shadowing algorithm. Through the parameter optimization 
procedure described in section \ref{sec:optim}, we find an optimal parameter 
that minimizes the observation error. This amounts to finding an $\epsilon$ such that 
the error between an observation orbit and a true orbit of $f(\cdot, \epsilon)$ is minimized.
Then, we can refine the background trajectory, using the shadowing perturbation (computed by the shadowing algorithm). The parameter optimization is then repeated 
with the refined trajectory as the new background. By repeating this combined 
state-parameter optimization procedure, we minimize both parts of the 
state estimation cost functional separately. At the end of this procedure, 
the analysis state is obtained by iteratively refining the background state using 
the shadowing perturbations at different parameters. 
We outline the state-parameter optimization procedure 
assuming we have observations at every timestep. 
The following steps are repeated, until the cost functional is 
less than a specified tolerance. Initially, 
the reference orbit $u_n$ for shadowing is set to the 
background orbit. The algorithm is the following:
\begin{enumerate}
    \item Run tangent shadowing (section \ref{sec:nilss}) with an objective function 
    \begin{align}
    \langle J\rangle(s) := \dfrac{1}{N}\sum_{n=0}^{N-1}(g^{\rm obs}_n - g_n)^2,
    \end{align} to obtain $d_s \langle J\rangle,$ and the sequence $\{v^{\rm sh}_n\}.$
    \item Update the parameter as: $s \longleftarrow s + \delta s$, where 
    $\delta s = \gamma \;d_s \langle J\rangle.$
    \item Update the initial condition for the next iteration 
    as $u_0 \longleftarrow u_0 + \delta s \; v^{\rm sh}_0$. Go to step 1. 
\end{enumerate}
In practice, the time $0$ corresponding to the start of the assimilation window, must be postponed roughly by $1/\lambda_1$ (Lyapunov time), to allow a spin-up time for $v^{\rm sh}_0$ to be accurately computed. 
\begin{figure}
    \centering
    \includegraphics[width=0.9\textwidth]{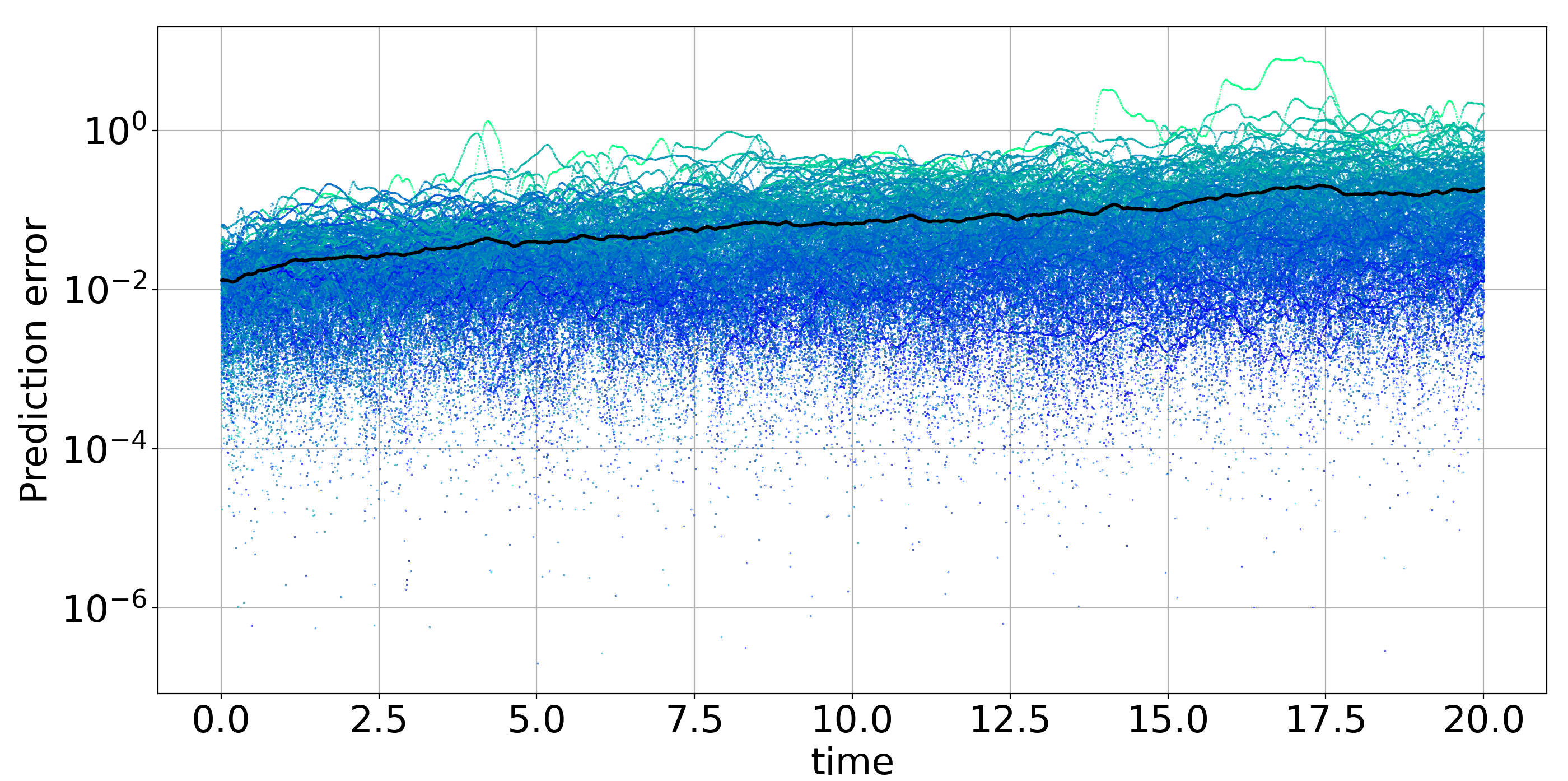}
    \caption{Relative error in the instantaneous acoustic energy between its predicted and observed values in the chaotic regime ($\beta = 7$). The maximum error increases from blue to green colors. The mean error (over the assimilation time window) is shown in black.}
    \label{fig:rijke_pred_error}
\end{figure}
The relaxation factor $\gamma$ is assumed to be a fixed constant, as in the parameter optimization procedure in section \ref{sec:optim}.
\subsection{Numerical results}
We present state estimation results on the time-delayed model, computed using the above algorithm. The algorithm is validated on the 
Lorenz'63 model, as shown in Figure \ref{fig:lorenz63_pred_error} in \ref{sec:appxLorenzDA}. 
A total of 180 experiments are performed on the time-delayed model, each with a different background trajectory of length 20 time units. The parameter $s$ that is updated is set to 
$\beta,$ with the reference value of 7.0. In each experiment, a background state is generated by perturbing each component of a reference state 
by a Gaussian random variable of variance 0.1. The original trajectory, started from the reference state,
is used to generate the observations (of the acoustic energy) $J_{\rm ac}^{\rm obs}$ at every timestep. The mean squared observational 
error is the objective function for 
tangent NILSS, $$ \langle J\rangle := \dfrac{1}{N}\sum_{n=0}^{N-1} \lvert J_{{\rm ac}_n}^{\rm obs} - J_{{\rm ac}_n}\rvert^2.$$ 
The step size for gradient descent (section \ref{sec:optim}) is $\gamma = 0.1$. We show in Figure \ref{fig:rijke_pred_error} the relative errors in the acoustic energy along the analysis orbit, which is the result of the algorithm after 200 gradient descent steps. The relative error is defined as 
$\lvert J^{\rm obs}_{{\rm ac}_n} - J_{{\rm ac}_n}\rvert/J^{\rm obs}_{{\rm ac}_n}$. Each colored line indicates a single experiment, with a total of 180 numerical 
experiments performed with different background states. The color of the line
indicates the maximum relative error observed in that experiment,
during the assimilation window of 20 units; the maximum errors 
increase from blue to green. The sample mean of the relative errors across all the 
experiments is shown in black. Note that the assimilation time is 4 times longer 
than $1/\lambda_1,$ the Lyapunov time. 
As shown in Figure \ref{fig:perts}, over 20 time units, we expect 
a small perturbation introduced in the initial condition, in almost any direction, to
grow by 3 orders of 
magnitude. However, as the results in Figure \ref{fig:rijke_pred_error} indicate, the relative error in $J_{\rm ac}$  has been restricted to within 
10 \% over this assimilation window, due to use of shadowing directions to iteratively refine the initial condition and the parameter, to match the 
observations. Improvements to this algorithm in order to reduce the errors, and further increase the predictability window will be studied in future work. 
One modification 
to the suggested algorithm, toward this goal, is to incorporate the observation error at every timestep, into the perturbation $x_n$ in 
the shadowing algorithm. 

\section{Conclusions}
\label{sec:conclusions}
Na\"ive applications of linear perturbation methods such as tangent/adjoint/Automatic Differentiation (AD)/finite-difference, cannot compute 
the derivatives of long-time averages in chaotic systems with respect to specified inputs. A 
recent method, known as the Non-Intrusive Least Squares Shadowing \cite{angxiu1}, computes these derivatives 
by numerical construction of tangent/adjoint \emph{shadowing} perturbations, which are infinitesimal 
perturbations that remain bounded for a long time duration. In this paper, we introduce AD into the tangent and 
adjoint NILSS. This nontrivial combination of algorithms is an enabler for the application of shadowing to complex dynamical systems. We demonstrate shadowing by computing sensitivities on a chaotic time-delayed model, which 
is a reduced-order model of a gas turbine combustor. We compute tangent and adjoint shadowing sensitivities of the 
ensemble averages of the acoustic energy and Rayleigh index with respect to the design parameters that control the heat release rate and the time delay.
Although the model is a reduced representation compared to CFD-based combustion models, it can be used to estimate the development of chaotic acoustic instabilities. 
First, we demonstrate an automatically-differentiated  procedure  for the minimization of the long time-averaged acoustic energy through heat-release parameter selection, which does not require a tangent solver~\cite{Huhn2020_jfm}. Secondly, we construct a  pseudo-orbit data assimilation scheme using the computed shadowing sensitivities in an optimization loop. We show that this scheme extends the predictability window by four Lyapunov times. Finally, we remark that the proposed algorithm, 
the shadowing-based optimization and data assimilation scheme are more generally applicable. The algorithms and the software  developed (available at 
\cite{nisha-energies}) may be used for other hyperbolic chaotic models with/without time delay. 

\section*{Funding}
Nisha Chandramoorthy and Qiqi Wang gratefully acknowledge the support of the Air Force Office of Scientific Research Grant No. FA8650-19-C-2207. Luca Magri gratefully acknowledges the support of the Royal Academy of Engineering Research Fellowships Scheme. \textbf{Conflict of Interest}: The authors declare that
they have no conflict of interest. 
\appendix

\section{Validation of sensitivity computation on the Lorenz'63 model}
\label{sec:appxLorenzSens}
\begin{figure}
	\includegraphics[width=0.5\textwidth]{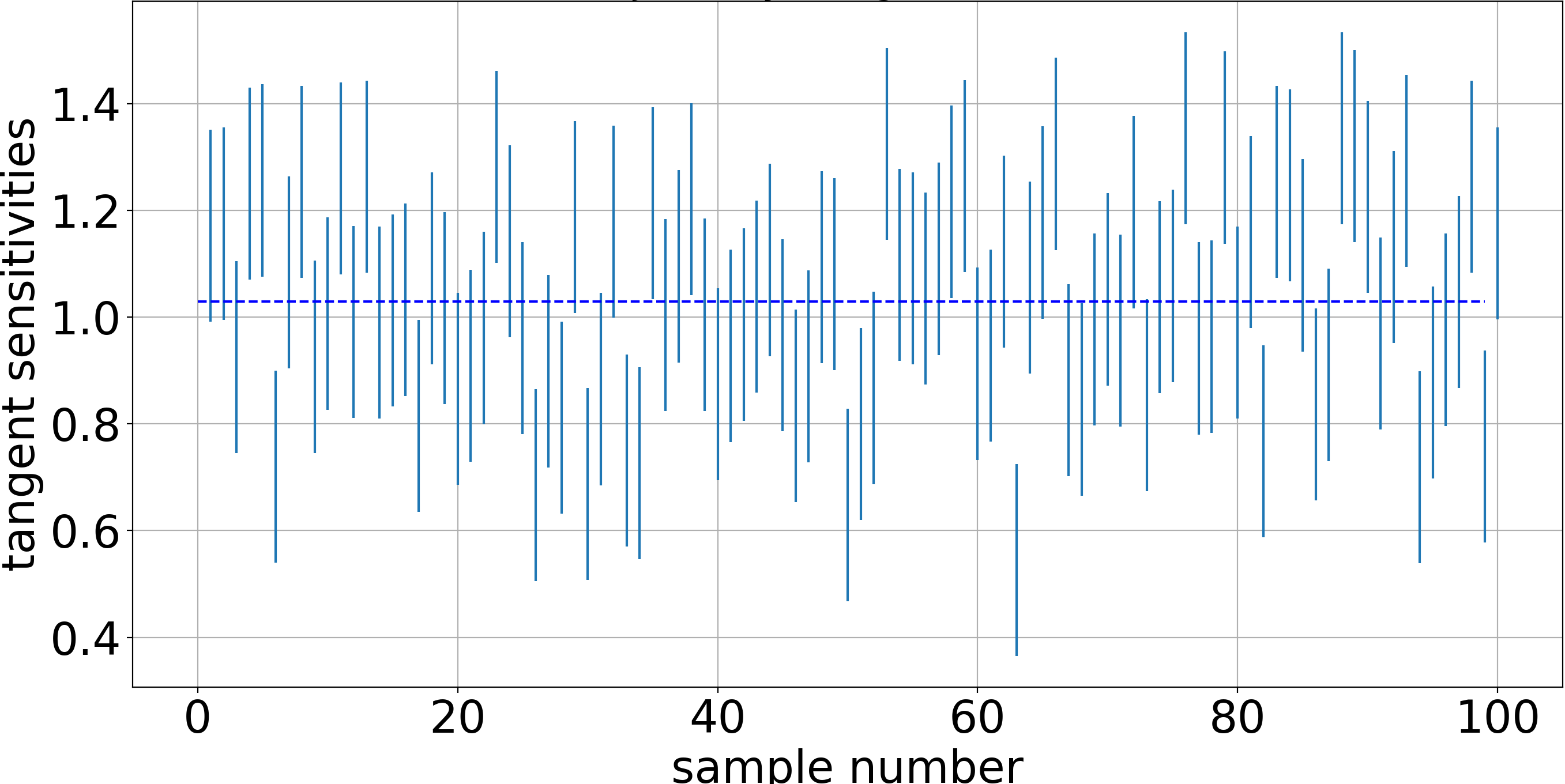}
	\includegraphics[width=0.5\textwidth]{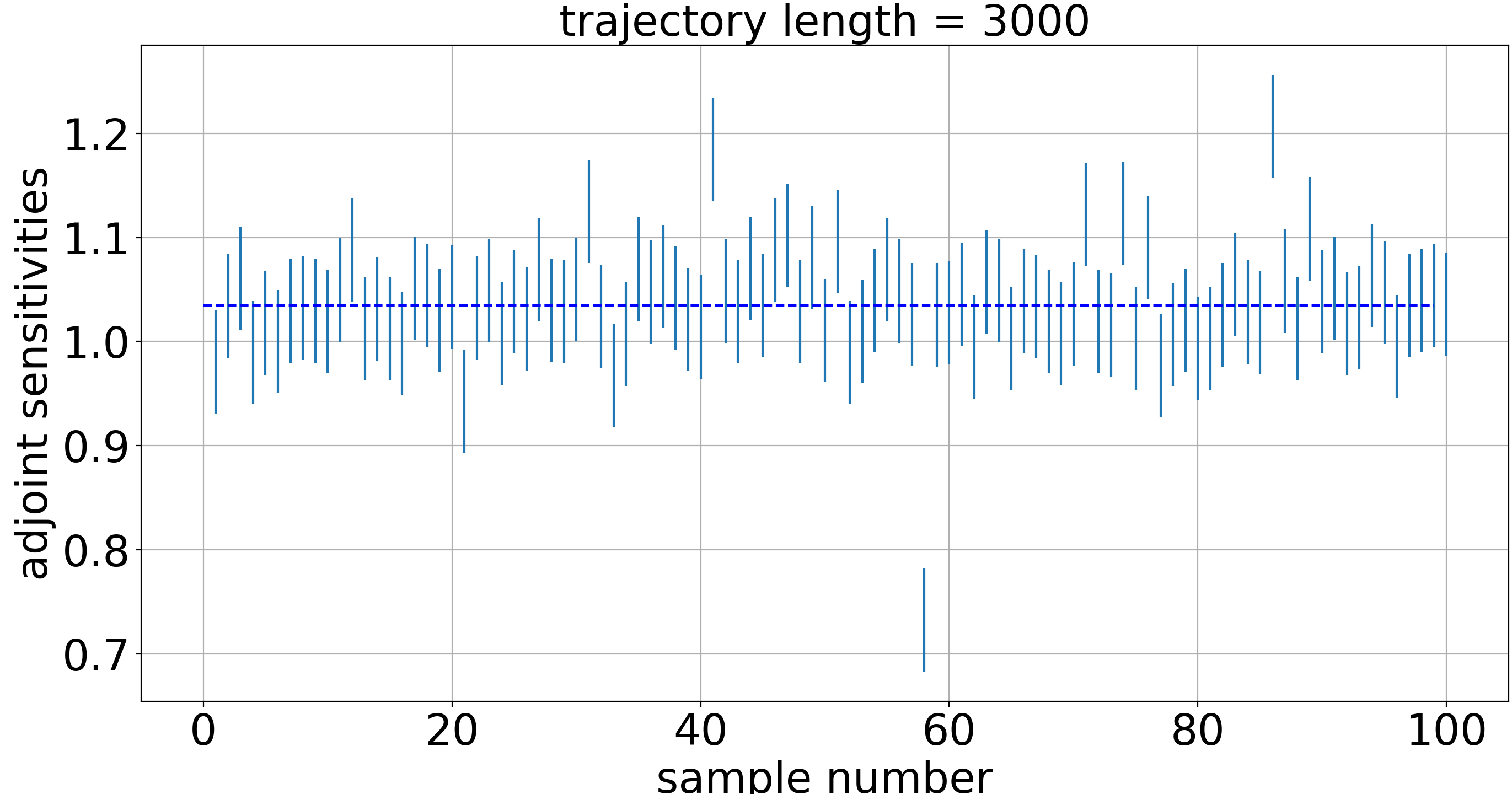}
	\caption{
	Sensitivities $d \langle z \rangle/ds$ computed for the Lorenz'63 model using the tangent algorithm (left) and the adjoint algorithm (right). Different trajectories of length 3000 or 15 timeunits are used to compute the sensitivities, shown as an errorbar of length one standard deviation.}
	\label{fig:lorenz63_sensitivities}
\end{figure}

We use the classical model of chaos, the Lorenz'63 system, for validation results in this paper. The Lorenz'63 model is the following
set of nonlinear ODEs that serves as a reduced-order model for fluid 
thermal convection between parallel plates maintained at a temperature difference \cite{lorenz}:
\begin{align}
		\dfrac{d}{dt} \begin{bmatrix}
				x\\
				y\\
				z
		\end{bmatrix}
		&= \begin{bmatrix}
				10(y-x) \\
				x(s - z) - y\\
				xy - (8/3) z
			\end{bmatrix}.
	\label{eqn:lorenz}
\end{align}
Here a phase point $u$ is represented using 3 coordinates as $u \equiv [x,y,z]^T \in \mathbb{R}^3$ and the map $f(\cdot, s)$ advances a state $u_n$ to $u_{n+1}$ by timestepping the ODE system (Eq. \ref{eqn:lorenz}). For the time integration, we use a Forward Euler scheme with a timestep of 0.005. The Lyapunov exponents of this system are about 0.9, 0, and -14.6. This is a partially hyperbolic system \cite{lorenz-partial-hyperbolicity} that has been shown to possess an SRB measure \cite{lorenz-CLT}. 
The objective function for validating shadowing is chosen to be the $z$ coordinate function. It is known that the
ergodic average $\langle z\rangle$ as a function of $s$, can be 
approximated as a straight line with a slope of about 0.96, over a range of values around the standard $s = 28$ \cite{Lea2000}. We use the discrete shadowing algorithms described in section \ref{sec:nilss} in both tangent and adjoint mode, to compute this derivative. Numerical results are shown 
in Figure \ref{fig:lorenz63_sensitivities}, in which 100 sample sensitivities are shown each over a time duration of 15 units, starting from different points on the Lorenz'63 attractor. On the left, we see the computed tangent shadowing sensitivities with $d_2 = 2$, and the on the right, the adjoint shadowing sensitivities, also computed with $d_u = 2$. From the numerical results, we see that the sample means of both sensitivities are within 10\% of the reference value of 0.96, thus validating both tangent and adjoint shadowing sensitivity codes. These tests can be run from \verb+tests/test_lorenz63.jl+ at \cite{nisha-energies}. Interestingly, adjoint shadowing appears to be better suited for this objective function-parameter pair since the variance of tangent sensitivities is 10 times larger than the adjoint.

\section{Validation of data assimilation scheme on the Lorenz'63 model}
\label{sec:appxLorenzDA}
\begin{figure}
    \centering
    \includegraphics[width=0.9\textwidth]{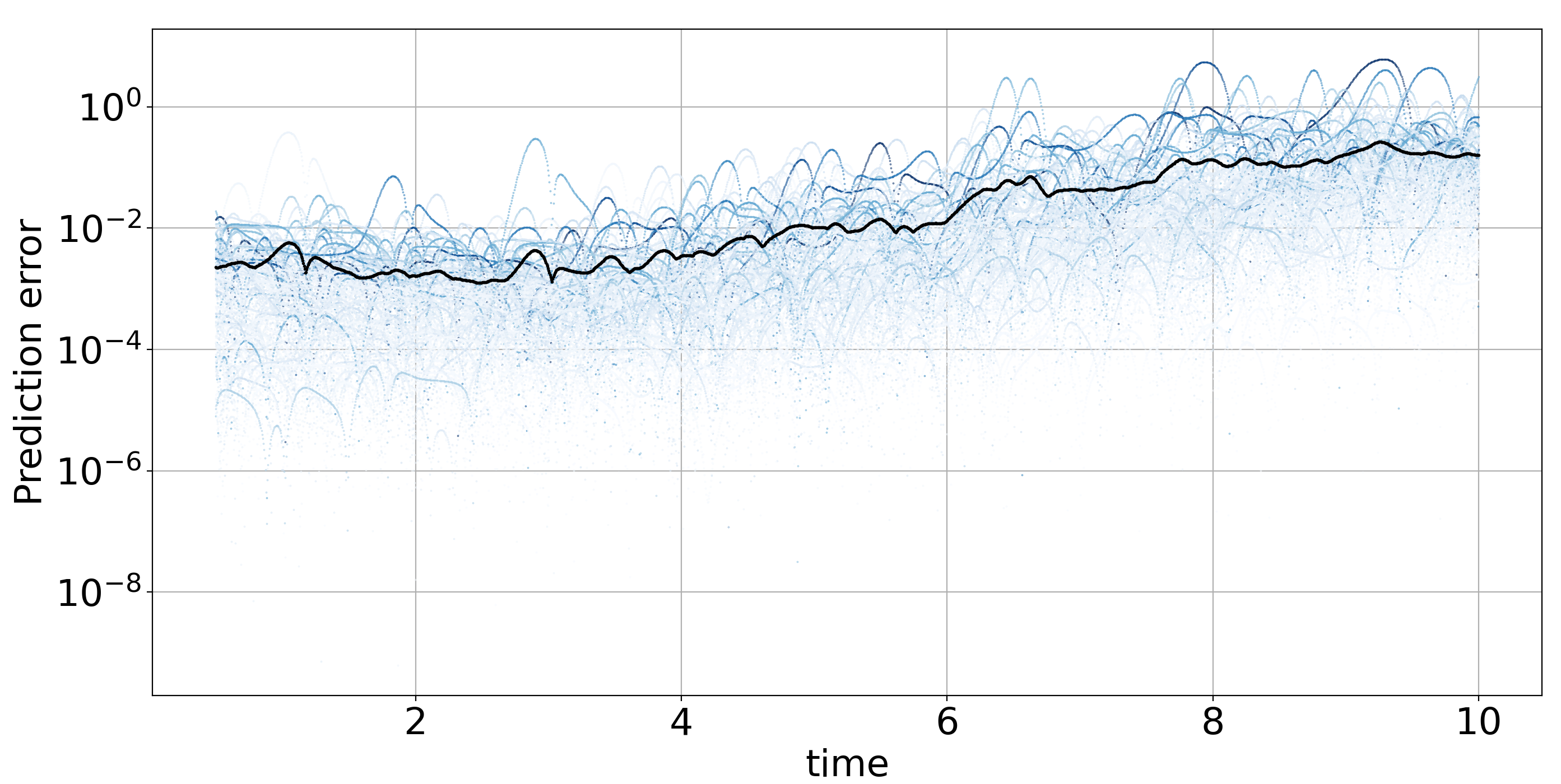}
    \caption{Relative error in the predicted state $z$ as a function of time for the Lorenz'63 system. The maximum prediction error increases from light to dark colors. The mean error across all experiments is shown in black.}
    \label{fig:lorenz63_pred_error}
\end{figure}
We apply the data assimilation scheme described in section \ref{sec:da} on the Lorenz'63 model. We perform a series of 100 numerical experiments each with a different background. We fix a reference trajectory of length 5000  steps (10 time units), and generate a background by perturbing  
the reference initial condition in the $\hat{z}$ direction by an additive Gaussian noise with a variance of 0.1. The reference trajectory is 
used to generate the observation trajectory $J^{\rm obs}_n := z^{\rm obs}_n.$ Tangent NILSS 
is run with an objective function,
\begin{align}
    \langle J\rangle(s) := \dfrac{1}{N}\sum_{n=0}^{N-1}(z^{\rm obs}_n - z_n)^2,
\end{align}
where $z_n$ is the $z$-coordinate of a trajectory $u_n.$
Using the derivative $(d\langle J\rangle/ds)(s)$ computed by shadowing,
the parameter $s$ is updated in a gradient descent algorithm with a 
constant step size of $\gamma = 0.1$. At the beginning of gradient descent, the trajectory is set to the background trajectory. At each gradient descent step, it is iteratively refined 
to a pseudo-orbit, as described in Step 3 in section \ref{sec:da}. 
The mean error is calculated by averaging the optimal error across the 100 experiments. The optimal errors 
, defined as relative observation errors after 200 gradient descent steps, across experiments are as shown as a function of time in 
Figure \ref{fig:lorenz63_pred_error}. The relative observation error at time $n$ is 
$|(z^{\rm obs}_n - z_n)|/z^{\rm obs}_n$. The colors of the lines in Figure \ref{fig:lorenz63_pred_error} are 
according to the maximum over $n$ of the relative observation error. As mentioned before, the only positive Lyapunov exponent of this system is known to be about 0.9. The results indicate a predictability within 10\% of the observation, 
on average, even up to 10 time units, which is about $10 \times (1/\lambda^1).$ With this validation on the Lorenz'63 model, we apply the same method for data assimilation in the Rijke tube model, as described in section \ref{sec:da}.

\section{AD shadowing and code for replication/extension}
\label{sec:appxAD}
As we describe in the main text, the AD version of shadowing is enabled by the introduction of AD to replace the tangent/adjoint solvers in tangent/adjoint shadowing. In this section, we give a pseudocode for this modification. The source 
code in Julia for generating the numerical results in this paper is available at \cite{nisha-energies}. This section also briefly describes how this code may be used with a different chaotic ODE/map.
We must mention that an alternative AD-version of the NILSS algorithm is also available at \cite{gitnisha}, in which the \verb+OpenAD+ AD software package \cite{openad,openAD1,paul-AD}, written in Fortran, is used. This latter code assumes that the primal solver is available from the user as a binary file.

In the AD version of tangent/adjoint shadowing, as mentioned before, we replace the needed tangent/adjoint solvers with AD. In order to do this, we define a function, say \verb+f+, whose return value will be differentiated in forward/reverse-mode, and we specify some input variable(s) to differentiate with respect to. The AD software returns the gradient of the return value with respect to the input which is a function of the arguments passed to \verb+f+. Let \verb+f+ be a function 
with the arguments $x,y$. The scalar or vector return value 
of this function is written as $\verb+f+(x,y)$. 
For instance, using the \verb+Zygote.jl+ Julia package, 
one would use the following syntax to obtain the gradient of $\verb+f+(x,y)$ with respect to $x$:
$\verb+Zygote.gradient+(x \verb+-> f+(x,y), x)$. If using
the \verb+OpenAD+ Fortran package instead, we declare the 
input $x$ and the output of \verb+f+ to be \emph{active} variables. If \verb+ret_f+ is the variable that stores the return value of $\verb+f+$, its derivative with respect to $x$ is stored in the variable \verb+ret_f%d+. In the AD version of tangent/adjoint shadowing, we only propose to modify the $n$-loop in section \ref{sec:nloop} by introducing AD to replace the tangent/adjoint solvers; the rest of the algorithm is retained as described in section \ref{sec:nilss}. The inputs to the AD shadowing algorithm are the sequences $\left\{u_n\right\}$ and $\left\{J_n\right\}$. We now give the modified $n$-loop in which we define the needed function, and input variables, for each AD invocation.
\begin{enumerate}
		\item Obtain $q_n^i$ from $q_{n-1}^i$ by applying AD. In 
				particular, 
			\begin{itemize}
					\item for tangent shadowing, 
						we differentiate the function 
						$f(u_{n-1} + \epsilon q_{n-1}^i,s)$.
					The input variable is $\epsilon$ and the output 
				is the return value of the function. The derivative 
							obtained from forward-mode AD is $q_n^i$, and this 
							must be carried out for $1\leq i\leq d_u$;
    				\item  for adjoint shadowing, 
						we differentiate the value 
							$f(u_{n'+1},s) \cdot q_{n-1}^i$, where $n' = N + 1 - n$.
							The input variable is $u_{n'+1}$. 
							The derivative 
							obtained from reverse-mode AD is $q_n^i$, and this 
							must be carried out for $1\leq i\leq d_u$.

			\end{itemize} 
    \item Let $Q_n$ be the matrix with the columns $q_n^i$, $1\leq i\leq d_u$. 
    QR-factorize $Q_n$ and set $Q_{n}$ to the 
    obtained ``Q''. Let the ``R'' from QR factorization be stored 
    as $R_{n}$. Thus, each $q_n^i$, $1\leq i \leq d_u$ is now a 
    unit vector.
    
    \item Obtain $v_n$ from $v_{n-1}$ by applying AD. In particular,
			\begin{itemize}
				\item for tangent shadowing, 
						we differentiate the function 
						$f(u_{n-1} + \epsilon v_{n-1},s + \epsilon)$.
					The input variable is $\epsilon$. The derivative 
							obtained from forward-mode AD is $v_n$;
    				\item  for adjoint shadowing, 
						we differentiate the value 
							$v_{n-1}\cdot f(u_{n'+1}, s) + 
							(1/N) J_{n'+1}$.
							The input variable is $u_{n'+1}$. 
							The derivative 
							obtained from reverse-mode AD is $v_n$.
			\end{itemize}
    
    \item Set $\pi_n := Q_n^T v_n$, which is a $d_u$-length 
    orthogonal projection row-vector of $v_n$ along $q_n^i.$
    
    \item Normalize $v_n$ by projecting out the components 
    along $q^i_n$, that is, set $v_n \to v_n - \pi_{qv} Q_n.$
    
    \item Go to step 1 with $n\to n + 1$, and repeat until 
    $n = N.$
\end{enumerate}
In the code at \cite{nisha-energies}, the Julia script \verb+tests/test_lorenz63.jl+ and \verb+tests/test_rijke.jl+ compute respectively, for the Lorenz'63 model (Eq. \ref{eqn:lorenz}) and the Rijke tube model (Eq. \ref{eqn:rijke}),
the tangent and adjoint shadowing sensitivities. In order to use the code 
at \cite{nisha-energies} with a different model, we can use one of these test files as a template. The Julia file describing the model equations, which must be included in the test file, can follow the existing examples in \verb+examples/lorenz63.jl+ or \verb+examples/rijke.jl+. The optimization and data assimilation routines described in section \ref{sec:optim} and \ref{sec:da} respectively, 
can also be applied to a new model. To do this, we include the file containing the 
model equations into the utilities i) \verb+utils/optimize.jl+ for parameter optimization and 
ii) \verb+utils/rijke_tangent_state_estimation.jl+ for the state estimation algorithm.  
\bibliographystyle{elsarticle-num} 
\bibliography{main,MyCollection}

\end{document}